\documentclass[12pt]{amsart}
\usepackage{amscd,amssymb}
\theoremstyle{plain}
\newtheorem{Rmk}[equation]{Remark}
\newtheorem{Thm}[equation]{Theorem}
\newtheorem{Ex}[equation]{Example}
\newtheorem{Cor}[equation]{Corollary}
\newtheorem{Prop}[equation]{Proposition}
\newtheorem{Lem}[equation]{Lemma}
\newtheorem{Def}[equation]{Definition}
\newtheorem{Fact}[equation]{Fact}
\numberwithin{equation}{section}
\newfont{\goth}{eufm10 at 12pt}

\def\bt{\begin{Thm}}
\def\et{\end{Thm}}
\def\br{\begin{Rmk}}
\def\er{\end{Rmk}}

\def\bc{\begin{Cor}}
\def\ec{\end{Cor}}
\def\bp{\begin{Prop}}
\def\ep{\end{Prop}}
\def\bl{\begin{Lem}}
\def\el{\end{Lem}}
\def\bd{\begin{Def}}
\def\ed{\end{Def}}
\def\bq{\begin{quotation}}
\def\eq{\end{quotation}}
\def\beq{\begin{equation}}
\def\eeq{\end{equation}}
\def\bea{\begin{eqnarray}}
\def\eea{\end{eqnarray}}
\def\bfa{\begin{Fact}}
\def\efa{\end{Fact}}
\def\bex{\begin{eqnarray*}}
\def\eex{\end{eqnarray*}}
\def\ba{\backslash}

\newcommand{\mc}[1]{{}}

\newcommand{\g}[1]{{\mbox{\goth #1}}}
\newcommand{\m}[1]{\mathbb{ #1}}
\renewcommand{\c}[1]{\mathcal{ #1}}

   \def\e{\`e}      
         
   \def\wt{\widetilde}
\def\al{\alpha}       \def\be{\beta}        \def\ga{\gamma}
\def\de{\delta}       \def\eps{\varepsilon}  
            
       \def\la{\lambda}      \def\rh{\rho}
\def\si{\sigma}                
\def\ph{\varphi}               
\def\om{\omega}       \def\Ga{\Gamma}       
             \def\Ph{\Phi}
                  \def\Om{\Omega}
\def\ba{\backslash}

\def\ra{\rightarrow}
\def\hb{\hfill$\Box$}

\def\noi{\noindent}

\def\vs{\vspace{1em}}
\def\D{\displaystyle}
\def\proof{\noi {\it Proof. }}

\def\q{{\m Q}}
\def\G{\Gamma}

\def\and{{\rm and}}

\newcommand{\op}{\operatorname}
\newcommand{\bg}{{\bf G}}
\newcommand{\bh}{{\bf H}}
\newcommand{\SL}{\operatorname{SL}}

\renewcommand{\H}{\operatorname{H}}

\newcommand{\vol}{\operatorname{vol}}
\renewcommand{\e}{\varepsilon}
\newcommand{\Vol}{\operatorname{vol}}

\newcommand{\z}{\m Z}
\newcommand{\n}{\mathbb N}
\begin{document}

\title[$S$-integral points]{ Effective equidistribution
of $S$-integral points on symmetric varieties}

\author{Yves Benoist and Hee Oh}

\address{DMA-ENS 45 rue d'Ulm Paris 75005}


\email{benoist@dma.ens.fr}

\address{Mathematics 253-37, Caltech, Pasadena, CA 91125
\newline current address: Mathematics, Box 1914,
Brown University, Providence, RI}

\email{heeoh@math.brown.edu}
\thanks{the second author is partially supported by NSF grant 0629322}

\begin{abstract}{
Let $K$ be a global field of characteristic not $2$.
Let ${\bf Z}={\bh}\ba {\bg}$ be a symmetric variety defined over $K$
and $S$ a finite set of places of $K$.
We obtain counting and equidistribution results
for the $S$-integral points of ${\bf Z}$.

Our results are effective when $K$ is a number field.
}\end{abstract}

\maketitle



\section{Introduction}
\label{secintro}

\subsection{General overview}
\label{secgeneral}
Consider a finite system of polynomial equations with integral
coefficients. Its set of solutions
defines an arithmetic variety ${\bf Z}\subset \mathbb C^d$
defined over $\mathbb Z$.
For a set $S$ of primes including the infinite prime
$\infty$, let $\mathbb Z_S$ denote the ring of $S$-integers of $\q$,
that is, the set of rational numbers whose denominators
are products of primes in $S$.
If $S=\{\infty\}$, $\mathbb Z_S$ is simply the
ring of integers $\mathbb Z$,
and if $S$ consists of all the primes, then $\mathbb Z_S$ is the
field $\mathbb Q$ of rational numbers.
For any subring $R$ of $\mathbb C$,
we denote by ${\bf Z}_R$ the set of
 points in ${\bf Z}$ with coordinates
in $R$.
One of the fundamental questions in number theory is
to understand the
properties of sets ${\bf Z}_{\mathbb Z_S}$.
In this paper, we obtain effective counting and equidistribution
results of the $S$-integral points, for $S$-finite,
in the case when ${\bf Z}$ is a symmetric variety.

The counting question in this set-up
has been completely
solved for the {\it integral} points via several different methods.
The first solution
 is due to Duke, Rudnick and Sarnak in 1993 \cite{DRS}
and their proof uses the theory of automorphic forms.
Almost at the same time, Eskin and McMullen
gave the second proof utilizing mixing properties of semisimple
real algebraic groups \cite{EM}.
 The third proof, due to Eskin, Mozes and Shah \cite{EMS}, is based on the
ergodic theory of flows on homogeneous spaces, more
precisely, Ratner's work on the unipotent flows.

The approach of \cite{EM} using mixing properties has several advantages
over the others in our viewpoint.
First it does not require the deep theory
of automorphic forms, avoiding technical difficulties in dealing with
the Eisenstein series as in \cite{DRS}. Secondly,
 although this was never addressed in \cite{EM},
in principle it also gives a rate of
 convergence which the ergodic method of \cite{EMS} does not give.
Thirdly the method can be extended to other global fields
of positive characteristic, which is again hard to be achieved
via the ergodic method.

For these reasons,
we develop the approach
of Eskin and McMullen \cite{EM} in this paper
in order to obtain {\it effective} results
 for the general $S$-integral points on symmetric varieties.

We use the mixing properties of $S$-algebraic semisimple groups,
with a rate of convergence.
 Implementing this in the counting problem,
a crucial technical ingredient is to verify
 certain geometric property, which was named the wavefront property
by \cite{EM}, for an $S$-algebraic symmetric variety.
We prove this using the
 polar decompositions for non-archimedean symmetric spaces obtained
 in \cite{BO} specifically for this purpose.
We emphasize that the wave front property is
precisely the reason that our proofs work in the setting of
an $S$-algebraic symmetric variety for $S$ finite.
This property  does not hold
 for a general homogeneous variety even over the reals.
In obtaining effective counting results for
the $S$-integral points of bounded height, we also use
the works of Denef on $p$-adic local zeta functions
(\cite{De1}, \cite{De2}) and
 of Jeanquartier on fiber integrations \cite{Je}.

We remark that the approach for counting via mixing
was initiated in 1970 by Margulis in his dissertation
on Anosov dynamical systems \cite{Ma1}.
Recently similar mixing properties
in an adelic setting  have
been used in the study of rational points of {\it group} varieties
(see \cite{COU}, \cite{GMO}, and \cite{Gu}).
We also mention that for the case of group varieties,
the effective counting result was obtained for integral
points in \cite{GN} and \cite{Mau}.
We refer to  \cite{Le}, \cite{GW}, \cite{EO},
\cite{EV}, \cite{GOS1}, \cite{EL},
\cite{MV} etc.,
for other types of counting and equidistribution results.

\subsection{Main results}
\label{secresult}
We now give a precise description of the main results of this paper.

Let $K$ be a global field of characteristic not $2$,
i.e. a finite extension
of $\m Q$ or of $\m F_q(t)$ where $q$ is an odd prime.
Let ${\bf Z}$
be a symmetric variety in a vector space ${\bf V}$ defined over $K$.
That is, there exist
a connected algebraic
almost $K$-simple group ${\bf G}$ defined over $K$,
a $K$-representation
$\rho:{\bf G}\ra {\bf GL}({\bf V})$
with finite kernel and
a non-zero point $z_0\in {\bf V}_K$
whose stabilizer ${\bf H}$ in ${\bf G}$ is a symmetric $K$-subgroup
of ${\bf G}$
such that ${\bf Z}=z_0{\bf G} $.
By a symmetric $K$-subgroup of ${\bf G}$, we mean a $K$-subgroup
whose identity component
 coincides with the identity component of
the group of fixed points ${\bf G}^\si$ for
an involution $\si$ of ${\bf G}$ defined over $K$.
We assume that the identity component
${\bf H}^0$ has no non-trivial $K$-character.

We fix a basis of the $K$-vector space
${\bf V}_K$ so that one can define,
for any subring $\c O$ of $K$,
the subsets  ${\bf V}_{\c O}\subset {\bf V}$,
 ${\bf Z}_{\c O}\subset {\bf Z}$ and ${\bf G}_{\c O}\subset {\bf G}$
of points with coefficients in $\c O$.
For each place $v$ of $K$,
denote by $K_v$ be the completion of $K$ with respect to
the absolute value $|\cdot|_v$.
We write ${\bf V}_v$, ${\bf Z}_v$ and
${\bf G}_v$ for ${\bf V}_{K_v}$,
${\bf Z}_{K_v}$ and ${\bf G}_{K_v}$, respectively.

Let $S$ be a finite set of places of $K$ containing
 all archimedean (sometimes called infinite)
 places with $\bg_v$ non-compact.
Note that if $\rm{char}\, K$ is positive, $K$
does not have any archimedean place.
We denote by $\c O_S$ the ring of $S$-integers of $K$, that is,
\begin{eqnarray*}
\label{eqnos}
\c O_S:=\{ k\in K\;\mid \;
|k|_v\leq 1 \;\; \text{for each finite }v\not\in S\}.
\end{eqnarray*}

For instance, if $K=\mathbb Q$,
we have $\c O_S=\mathbb Z_S$.
We set
$ {\bf Z}_S=\prod_{v\in S}{\bf Z}_{v}$
and similarly ${\bf G}_S$ and ${\bf H}_S$.

Note that the sets ${\bf Z}_{\c O_S}$,
${\bf G}_{\c O_S}$ and ${\bf H}_{\c O_S}$ are discrete subsets
of ${\bf Z}_S$, ${\bf G}_S$, and ${\bf H}_S$ respectively,
 via the diagonal embeddings.

By a theorem of  Borel and Harish-Chandra
in characteristic $0$ and of Behr and Harder
in positive characteristic
(see Theorem I.3.2.4 in \cite{Ma2}),
the subgroups ${\bf G}_{{\c O}_S}$ and
${\bf H}_{{\c O}_S}$ are lattices in
${\bf G}_S$ and ${\bf H}_S$ respectively.
Again, by a theorem of  Borel, Harish-Chandra, Behr and Harder,
the group
${\bf G}_{\c O_S}$ has only finitely many orbits
in ${\bf Z}_{{\c O}_S}$
(see Theorem 10 in \cite{Go}).

Hence our counting and equidistribution question of $S$-integral points
of ${\bf Z}$ reduces to counting and equidistribution of
points in a single ${\bf G}_{\c O_S}$-orbit, say, for instance,
in $z_0{\bf G}_{\c O_S}$.
 Set
$$Z_S:=z_0{\bf G}_S=\prod_{v\in S}z_0 {\bf G}_v  \; $$
and let $\Ga_S$ be a subgroup of finite index in ${\bf G}_{\c O_S}$.
Let $\mu_{X_S}$ be a ${\bf G}_S$-invariant measure on
$X_S:=\Ga_S\backslash {\bf G}_S$ and $\mu_{Y_S}$
an ${\bf H}_S$-invariant measure on
$Y_S:=( \Gamma_S\cap {\bf H}_S)\backslash {\bf H}_S$.
For each $v\in S$,
we choose an invariant measure
 $\mu_{Z_v}$ on $Z_v:=z_0{\bf G}_v$  so that
for $\mu_{Z_S}:=\prod_{v\in S}\mu_{Z_v}$,
we have $\mu_{X_S}=\mu_{Z_S}\mu_{Y_S}$ locally.
For a subset $S_0\subset S$,
we set $\mu_{Z_{S_0}}=\prod_{v\in S_0}\mu_{Z_v}$.

For a Borel subset $B$ of $Z_S$,
we set
$$\vol(B):=
 \tfrac{\mu_{Y_S}(Y_S)}{\mu_{X_S}(X_S)}\;\mu_{Z_S} (B) .$$

We assume that ${\bf G}_S$ is non-compact; otherwise
${\bf Z}_{\c O_S}$ is finite.
By considering a finite covering of $\bg$ by its simply connected
cover, we may also assume that
${\bf G}$ is simply connected without loss of generality.

Before stating our main result, we
summarize our set-up:

\bq
$K$ is a global field such that char$(K)\!\neq\! 2$,
 ${\bf Z}\simeq {\bf H}\backslash {\bf G}$ is a symmetric variety
in a vector space ${\bf V}$ defined over $K$
where
${\bf G}$
is an almost $K$-simple simply-connected
$K$-group acting on ${\bf V}$
such that the identity component
${\bf H}^0$ has no non-trivial $K$-character, and
$S$ is a finite set of places of
$K$ containing all the infinite places $v$ with $\bg_v$ non-compact
and satisfying that ${\bf G}_S$ is non-compact.
\eq

\subsubsection*{{\bf Counting $S$-integral points}}
We first state our counting results. We refer
to Definition \ref{defwellrounded} for the notion of
a well-rounded sequence  of subsets $B_n$ in $Z_S$.
Roughly speaking, this means that for all small
$\e>0$, the boundaries of $B_n$ can be approximated
by neighborhoods whose volume is of $\e$-order compared to
the volume of $B_n$ uniformly.
\bt
\label{thintro}
For any well-rounded sequence  of
subsets $B_n$ of $Z_S$
with volume tending to infinity,
we have
$$\# ( z_0\G_S \cap B_{n})
\sim \vol(B_n) \quad\text{as $n\to\infty$}. $$
\et

As a corollary,
we obtain that the number of
 $S$-integral points of size less than $T$ is given by the
volume of the corresponding ball in $Z_S$.
A most natural way to measure the size of an $S$-integral point
is given by a height function $\H_S$.
For $z\in {\bf Z}(\c O_S)$, it is simply
$$\H_S(z):=\prod_{v\in S} \|z\|_v$$
where the $\|\cdot \|_v$ are norms on ${\bf V}_{K_v}$ which are euclidean
when $v$ is infinite and which are max norms when $v$ is finite.
This height function $\H_S$ naturally extends to $Z_S$.

\bc
\label{corintro}
As $T\to \infty$,
$$\#\{z\in z_0{\G_S}: \H_S(z) <T\} \sim \Vol (B_S(T)).$$
where $B_S(T):=\{z\in { Z}_S: \H_S(z) <T\}$.
\ec

When $K$ is a number field,
Theorem \ref{thintro} is proved with
a rate of convergence (see Theorem \ref{thspeed}).
For instance, we get:

\begin{Thm} \label{ec1}
Let $K$ be a number field.
There exists $\delta>0$ such that as $T\to\infty$
$$\#\{z\in z_0\Gamma_S: \H_S(z) <T\} =\Vol (B_S(T)) (1+
O(T^{-\delta})). $$
\end{Thm}

We will see (Remark \ref{remvt}) that there exist $a\in \m Q_{>0}$, $b\in \m Z_{\geq 0}$
and $c_1,c_2>0$ such that for $T$ large,
$$
c_1 T^a\log (T)^b
\leq \Vol (B_S(T)) \leq
c_2 T^a \log (T)^b\; .
$$
In general, one cannot choose $c_1=c_2$.

The rate of convergence in Theorem \ref{ec1}
is new even for integral points in the generality of symmetric
varieties. In this case, as $T\to\infty$, $\Vol (B_S(T)) \sim
c\, T^a \log (T)^b$, for some $c>0$.

Our proof of Theorem \ref{ec1}
uses  Denef's result
on local zeta functions
which is not available in positive characteristic.
This explains our hypothesis on the characteristic of $K$.

\subsubsection*{{\bf  Equidistribution of $S$-integral points}}
To motivate, consider the case when $K=\q$
and suppose that ${\bf Z}_{\mathbb Z [p^{-1}]}$,
the set of rational points in ${\bf Z}$ with denominators only power of
$p$, is a dense subset
in ${\bf Z}_{\mathbb R}$, which is often the case.
A natural question is when
the sequence of subsets in ${\bf Z}_{\mathbb Z [p^{-1}]}$ consisting
of elements of denominator precisely $p^n$ is equidistributed
as $n\to \infty$.
That is, for two compact subsets $\Omega_1, \Omega_2$ of
${\bf Z}_{\mathbb R}$, as $n\to\infty$,
$$\frac{\{x\in {\bf Z}_{\mathbb Q}\cap \Omega_1: p^n x\in
 {\bf V}_{\mathbb Z}, \; p\nmid p^n x\}}
 {\{x\in {\bf Z}_{\mathbb Q}\cap \Omega_2:
p^n x\in
 {\bf V}_{\mathbb Z},\; p\nmid p^n x\}}
 \sim \frac{\Vol(\Omega_1)}{\Vol(\Omega_2)}?$$

Once we note that $p^n x\in V_{\mathbb Z}$ is equivalent to
the condition that the $p$-adic maximum norm of $x$
is at most $p^n$,
the above question can be rephrased as the question of
equidistribution on ${\bf Z}_{\mathbb R}$ of the sets
 $\{z\in {\bf Z}_{\mathbb Z[p^{-1}]}: \|z\|_p= p^n\}$.

We answer this question in greater generalities:
\bt
\label{thoeq}
 Let $S=S_0\sqcup S_1$ be a partition of $S$.
For any
well-rounded
sequence $B_n$ of subsets of $Z_{S_1}$
with volume tending to infinity, and
for any compact subset $\Omega\subset Z_{S_0}$
of positive measure and
of boundary measure $0$,
we have
$$\#\, z_0\Gamma_S\cap (\Omega\times B_n)\sim
 \tfrac{\mu_{Y_S}(Y_S)}{\mu_{X_S}(X_S)}
\; \mu_{Z_{S_0}} (\Omega)\, \mu_{Z_{S_1}}(B_n)
\quad\text{as $n\to\infty$}.$$
\et

Note that the special case discussed prior to Theorem \ref{thoeq}
corresponds to $K=\q$, $S_0=\{\infty\}$, $S_1=\{ p\}$, and
$B_n=\{z\in {\bf Z}_{\mathbb Q_p}: \|z\|_p= p^n\}$.

Note that in all the above
theorems, we may replace $z_0\G_S$ by
$Z_{\c O_S}:=Z_S\cap {\bf Z}_{\c O_S}$.
as long as we renormalize the volume form
so that
 the volume of a subset $E\subset Z_{S}$ is given by
\begin{eqnarray}
\label{eqnvolt}
\wt{\vol }(E)
&=&
\sum \tfrac{\mu_{Y_S}(Y_S)}{\mu_{X_S}(X_S)}\;\mu_{Z_S}(E)
\end{eqnarray}
where we sum the contributions from
each $\Ga_S$-orbit in
$Z_{\c O_S}$.
Hence we obtain:
\begin{Cor}
\label{coroeq}
Assume $S$ has at least two places.
\begin{enumerate} \item For any finite $v\in S$, the sets
 $Z(T):=\{z\in Z_{\c O_S}: \|z\|_v=T\}$
become equidistributed in ${ Z}_{S-\{v\}}$
as $T\to\infty$, subject to the condition
 $Z(T)\ne \emptyset$.
\item For an infinite $v\in S$,
the sets $Z_T:=\{z\in Z_{\c O_S}: \|z\|_v\leq T\}$
become equidistributed in ${ Z}_{S-\{v\}}$
as $T\to\infty$, provided ${ Z}_{v}$ is non-compact.
\end{enumerate}
\end{Cor}

Again, when $K$ is a number field,
Theorem \ref{thoeq} and Corollary \ref{coroeq} are proved with
a rate of convergence (see Corollary \ref{corspeed} and
Proposition \ref{prohbeffwell}). For instance, we obtain:

\begin{Thm} \label{thoeqef}
Let $K$ be a number field
and $S=S_\infty\sqcup S_f$ be the partition of $S$
into infinite and finite places. We assume that $\bg_{S_f}$ is non compact.
Set $\be_T:=\{ z\in Z_{S_f}\; :\; \H_{S_f}(z)<T\}$.
Then there exist $\delta>0$ such that
for any  compact subset $\Om$
of $Z_{S_\infty}$ with piecewise smooth boundary,
$$\#( z_0\Gamma_S\cap(\Om \times \be_T)) =
 \tfrac{\mu_{Y_S}(Y_S)}{\mu_{X_S}(X_S)}\; w_T\; \mu_{Z_\infty}(\Om) (1+
O(T^{-\delta}))\quad\text{ as $T\to\infty$} $$
where $w_T:=\mu_{Z_{S_f}}(\be_T)$.
\end{Thm}

\subsubsection*{{\bf Equidistribution of
 translates of $\bh_S$-orbits}}
Set  $X_S={\bf G}_{\c O_S}\ba {\bf G}_S$ and
$Y_S= {\bf H}_{\c O_S}\backslash {\bf H}_S$.
Let $\mu_{X_S}$ and $\mu_{Y_S}$ be invariant probability measures on
$X_S$ and $Y_S$ respectively.
The following theorem is a crucial tool in proving
Theorem \ref{thintro}. It states that
the translates $Y_Sg$
is equidistributed
in $X_S$
as $g$ leaves compact subsets of ${\bf H}_S\ba {\bf G}_S$.
\begin{Thm}
\label{eqqq}
For any $\psi\in  C_c(X_S)$,
$$\int_{Y_S}\psi(yg)d\mu_{Y_S}(y)\to \int_{X_S} \psi \;d\mu_{X_S}
\quad\text{as $g$ tends to infinity in ${\bf H}_{S}\ba {\bf G}_S$} .$$
\end{Thm}
The case when $K=\mathbb Q$ and $S=\{\infty\}$, Theorem \ref{eqqq}
 was proved in \cite{EM}, \cite{DRS},
\cite{EMS} and \cite{Sh}. In the case when $\bh_S$ is semisimple
and non-compact, it is recently
proved in \cite{GO}, by extending theorems of Mozes-Shah \cite{MS} and
Dani-Margulis \cite{DM} in $S$-algebraic settings.
None of the above papers address the rate issues, while our proof
gives effective version
in the case when $\text{char}(K)=0$: a
smooth function on $X_S$
is a function which is smooth for each infinite place in $S$ and
invariant under a compact open subgroup of ${\bf G}_{v}$
for each finite place $v\in S$.
The following effective version of theorem \ref{eqqq}
is a crucial tool in proving Theorem \ref{ec1} as well
as other effective results in this paper.
\begin{Thm}\label{transff}
For $K$ number field,
there exists $\kappa >0$ such that,
for any smooth function $\psi$ on $X_S$ with compact support,
there exists $c=c_\psi >0$ such that
 $$\left|\int_{Y_S} \psi (yg)\;
 d\mu_{Y_S}(y)-\int_{X_S}  \psi \;d\mu_{X_S}\right|
\le c\,\H_S( z_0g)^{-\kappa} \quad\text{ for all $g\in G_S$.}$$
 \end{Thm}

\subsubsection*{{\bf Examples}}
Let $K=\m Q$ and consider the following pairs $({\bf V}, f)$:
\begin{itemize}\label{invex}

\item[(A)] ${\bf V}$: the affine $n$-space with $n\ge 3$
and  $f$: an integral quadratic form
of $n$-variables. If $n=3$, we assume that $f$ does not represent $0$
over $\q$.
\item[(B)] ${\bf V}$: the space of symmetric $n\times n$ matrices
with $n\ge 3$ and $f=\pm \op{det}$.

\item[(C)] ${\bf V}$: the space of skew-symmetric $2n\times 2n$-matrices
with $n\ge 2$ and $f=\pm\op{pffaf}=\pm\sqrt{\op{det}}$.
\end{itemize}

For a positive integer $m$, define
$${\bf V}_m:=\{x\in {\bf V}: f(x)=m\} .$$
Consider the radial projection $\pi: {\bf V}_m \to {\bf V}_1$
given by $x\mapsto m^{1/d}x$ where $d$ is the degree of $f$.
Let ${\bf V}(\mathbb Z)^{\text{prim}}$
be the set of primitive integral vectors in ${\bf V}$.
For a finite set $S$ of primes of $\m Q$ containing the infinite prime
$\infty$, we denote by
 $\langle S \rangle \subset \m Q^*$
the multiplicative semigroup generated by the finite primes in $S$.

We give a partial answer to the following
Linnik problem (see \cite{Li}, \cite{Sa}, \cite{EO}, \cite{Oh2},
\cite{MV}):

\begin{Cor}
\label{corinv1}
Fix $S$
and $({\bf V}, f)$ as above.
Then
there exist constants $\delta>0$ and
 $\omega_m$,
such that
 for any non-empty compact subset $\Omega\subset {\bf V}_1(\m R)$
with piecewise smooth boundary,
we have
$$\# \,\Omega\cap \pi({\bf V}_m(\mathbb Z)^{\text{prim}})=
\omega_m \, \Vol(\Omega)
(1+O(m^{-\delta})) $$ as $m
\to \infty$ in  $\langle S\rangle$, subject to
${\bf V}_m(\z)^{\text{prim}}\ne \emptyset$.
\end{Cor}

A special case of (A) gives an effective equidistribution
for $\{x\in \z^3: f(x)=m\}$ with $f=x^2+y^2+z^2$
or $f=x^2+y^2-3z^2$, hence giving a different
proof of partial cases (because
of the restriction on $m$) of theorems of Iwaniec
\cite{Iw} and Duke \cite{Du}.
Note that a special case of (B) gives an effective equidistribution
for the {\it positive definite} integral matrices of given determinant.
These cases are of special interest
since the corresponding symmetric group $\bh$ is either
compact over the reals or a torus.
When ${\bh}$ is semisimple without compact factors over the reals,
Corollary \ref{corinv1} in its non-effective form, but
with no restriction on $m$, is obtained in \cite{EO}
using Ratner's work on the theory of unipotent flows.

\begin{Cor}
\label{corinv2}
Keep the same assumption as in Corollary \ref{corinv1}
and set $m_S$ to be the product of the finite $p\in S$.
For the case (A),
we further suppose
that $f$ represents $0$ over $\m Q_p$ for at least one $p\in S$.
Then
there exists $\delta>0$ such that
 $$\#\{x\in {\bf V}(\z)^{\text{prim}}:
 \|x\|_\infty<T,\;\; f(x)\in \langle S\rangle \}
 =v_T(1+O(T^{-\delta}))$$
where the asymptotic $v_T$ is given by the following
 sum over the divisors $m$ of
$m_S^{d-1}$:
$$\textstyle v_T=
\sum_{m|m_S^{d-1}}\wt{\vol}(\{x\in ({\bf V}_m)_S : \H_S(x)<T\})\;
 $$
\end{Cor}

To prove corollaries \ref{corinv1} and \ref{corinv2},
we will apply the effective versions of Theorems
\ref{thintro} and \ref{thoeq} to ${\bf V}_1(\mathbb Z_S)$.
We list more examples in section \ref{last}.

\subsection{Guideline}
\label{secplan}
We tried to help the reader in writing ''twice'' the proofs:
In the first half of this paper we concentrate
on the main term in the counting and equidistribution statements.
In the second half, we follow the same strategy
but develop more technical tools
to obtain the effective versions of these statements,
i.e. to control  the error terms.

In section \ref{secmixing}, we recall the decay of matrix coefficients
for semisimple groups $G$ and its application on a homogeneous
space $\Ga\ba G$ of finite volume.
In section \ref{secwavefront}, we show the wavefront property
for symmetric spaces $H\ba G$ over local fields and their products.
In section \ref{sechorbits}, we explain how mixing and wavefront
properties imply the equidistribution properties
of translates of $H$-orbits in $\Ga\ba G$ given in Theorem
\ref{eqqq}.
In sections \ref{secsumintegral} and \ref{seccounting},
we explain how these equidistribution properties
for translates of $H$-orbits in $\Ga\ba G$
allow us to compare
for {\it well-rounded} sequences of functions on $H\ba G$
each sum over a $\Ga$-orbit
with the integral on $H\ba G$.
In section \ref{secapplication}, we give examples of
well-rounded sequences and give proofs of Theorem \ref{thintro},
Corollary \ref{corintro}, Theorem
\ref{thoeq} and Corollary \ref{coroeq}.

Starting from section \ref{seceffmix},
we prove  the effective results listed in the introduction.
Theorem \ref{transff} is proved
 in section \ref{eff}, and  Theorems \ref{ec1}, \ref{thoeqef},
Corollaries \ref{corinv1} and \ref{corinv2},
among other effective applications, are
proved in section \ref{seceffapp}.
We list more concrete examples in section \ref{secexample}.

In the appendix \ref{secapp}, we give some general estimates
for the volume of balls in the orbits of algebraic groups
both over the real and p-adic numbers.

We remark that in the whole paper the assumption
of $\bh$ symmetric is used only to
obtain the (effective) wave front property for $\bh_S\ba \bg_S$.
The methods and the arguments in this paper work equally well for any
$K$-subgroup $\bh$ with no non-trivial characters
satisfying the wave front property.
\vs

\noindent{\bf Acknowledgment}
The  authors would like to thank
Alex Gorodnik for
helpful conversations.

\section{The mixing property}
\label{secmixing}

We first recall the Howe-Moore property also called {\it decay of
matrix coefficients}.
\bd
\label{defhowemoore}
A locally compact group $G$ is said to have the Howe-Moore
property if, for every unitary representation
$(\c H,\pi)$ of $G$ containing no non-zero vectors
invariant by a normal
non-compact subgroup, we have for all $v$, $w\in \c H$,
$$
\lim_{g\ra\infty}\langle \! \pi(g)v, w\! \rangle =0\; .
$$
\ed

This Howe-Moore property is related to the following mixing property.

Let $G$ be a (unimodular) locally compact group and $\Ga$ a lattice in $G$,
i.e. a discrete subgroup of finite covolume.
Let $\mu_X$ be a $G$-invariant measure on
$X:=\Gamma\ba G$. The group $G$ acts on $X$ by right-translations.
\bd
\label{defmixing}
The action of $G$ on $X$ is said to be {\it mixing}
 if for all $\al$ and $\be\in L^2(X)$
$$
\lim_{g\ra\infty}\int_X\al (gx) \be (x) d\mu_X(x) = \mu_X(X)
\int_X \al\; d\mu_X\;  \int_X \be\;d\mu_X .
$$
\ed

The relation between these two definitions is given by the
following straightforward proposition.

\bd
\label{defirreducible}
A lattice $\Ga$ in a locally compact group
$G$ is called ${\it irreducible}$ if
for any non-compact normal subgroup $G'$ of $G$, the subgroup $\Ga G'$
is dense in $G$.
\ed

Note that this definition is slightly
stronger than the usual definition since
it excludes lattices contained in a
proper subgroup of $G$.

\bp
\label{promixing}
Let $G$ be a locally compact group satisfying the Howe-Moore
property and $\Ga$ an irreducible lattice in $G$.
Then the action of $G$ on $\Ga\ba G$ is mixing.
\ep

\proof This is well-known. One may assume that $\al$ and $\be$
belong to $\c H:=L^2_0(X)$ of square-integrable
 functions with zero integral. The $G$-action by right-translations on
$\c H$ via $(\pi(g)f)(x)=f(xg)$ is a unitary representation of $G$.
The irreducibility hypothesis on $\Ga$ implies precisely
that $\c H$ does not contain
any non-zero vector invariant by a normal
non-compact subgroup of $G$.
Hence, by Definition \ref{defhowemoore},
the matrix coefficients $\langle \! \pi(g)\al,\be\! \rangle$ converge to $0$
as $g$ tends to infinity.
\hb\vs

The main example is due to Howe-Moore.

\bt
\label{prohowemoore}
For $i=1,..,m$, let $k_i$ be a local field and $G_i$  the group
of $k_i$-points of a connected semisimple $k_i$-group.
Then the product
$G:=\prod_{i=1}^m G_i$ has the Howe-Moore property.
\et

In this paper, ``local field'' means ``locally compact field'',
i.e. a completion of a global field, or, equivalently,
a finite extension of $\m R$, $\m Q_p$ or $\m F_p((t))$.

\proof See, for instance,
Proposition II.2.3 of \cite{Ma2} or
\cite{Be}.
\hb

\section{The wavefront property}
\label{secwavefront}

The wavefront property was introduced by Eskin and McMullen
for real symmetric spaces \cite{EM}.
Let $G$ be a locally compact group and $H$ a closed subgroup of
$G$.

\bd \label{defwavefront}
The group $G$ has the wavefront property in $H\ba G$ if
there exists a Borel subset $F\subset G$ such that $G=HF$
and, for every neighborhood $U$ of $e$ in $G$,
there exists a neighborhood $V$ of $e$ in $G$ such that
$$HVg\subset HgU
\quad\text {for all $g\in F$} .$$
\ed

This property means roughly that
the $g$-translate of a small neighborhood of
the base point $z_0:=[H]$ in $H\ba G$ remains
near $z_0g$ uniformly over $g\in F$.

This section is devoted to proving the following:

\bp
\label{prowavefront}
Let $k$ be a local field of characteristic not $2$,
${\bf G}$ a connected semisimple $k$-group, $\si$
a $k$-involution of ${\bf G}$, $G={\bf G}_k$ and $H$
a closed subgroup of finite index in the group
$G^\si$ of $\si$-fixed points.

Then the group $G$ has the wavefront
property on $H\ba G$.
\ep

To prove the above proposition, we need the following two lemmas.
A  $k$-torus ${\bf S}$ of ${\bf G}$
is said to be
$(k,\si)$-split if it is $k$-split
and if $\si(g)=g^{-1}$ for all $g\in {\bf S}$.
By a theorem of Helminck and
Wang \cite{HW}, there are only finitely many $H$-conjugacy classes
of maximal $(k, \sigma)$-split tori of ${\bf G}$.
Choose a set $\{{\bf A}_i:1\le i\le m\}$ of representatives of
$H$-conjugacy class of maximal $(k,\si)$-split tori
of ${\bf G}$ and set $A=\cup_{i=1}^m {\bf A}_i(k)$.

The following lemma was proved in \cite{EM}
for $k=\mathbb R$,
in \cite{BO} for all local fields of characteristic not $2$
(and independently in \cite{DS} when the residual characteristic
is not $2$).

\bl [{\bf Polar decomposition of symmetric spaces}] \label{lemkah}
There exists a compact subset $K$ of $G$ such that
$$
G=HAK\; .
$$
\el

The second lemma we need is based on the work of Helminck and Wang.

Let ${\bf A}$ be a maximal $(k,\si)$-split torus of ${\bf G}$
and ${\bf L}$ the centralizer of ${\bf A}$ in ${\bf G}$.
The set of  roots $\Ph=\Ph({\bf G},{\bf A})$
for the action of ${\bf A}$
on the Lie algebra of ${\bf G}$ is a root system.
For every positive root system $\Ph^+\subset \Ph$,
let ${\bf N}$ (resp. ${\bf N}^-$)
be the unipotent subgroup of ${\bf G}$ generated by
the root groups ${\bf U}_\al$
(resp. ${\bf U}_{-\al}$), for $\al\in \Ph^+$,
let ${\bf P}:={\bf L}{\bf N}$
(resp. ${\bf P}^-:={\bf L}{\bf N}^-$)
and ${\bf A}_k^+$ the Weyl Chamber~:
$$
{\bf A}_k^+:=\{ a\in {\bf A}_k \;|\;
|\al(a)|\leq 1\; ,\, \mbox{for all}\; \al\in \Ph^+\}\; .
$$
When $\Ph^+$ vary,
the Weyl chambers form a finite covering of ${\bf A}_k$.
Since ${\bf P}^-=\si({\bf P})$,
the parabolic $k$-subgroups ${\bf P}$ are $\si$-split,
i.e., the product ${\bf H}{\bf P}$ is open in ${\bf G}$
\cite[Prop. 4.6 and 13.4]{HW}.
Conversely, any
minimal $\si$-split parabolic $k$-subgroups of
${\bf G}$ containing ${\bf A}$ can be constructed in this way
for a suitable choice of $\Ph^+$.

\bl
\label{lemawaw}
\begin{enumerate}
\item The multiplication map
$m: {\bf H}_k\times{\bf P}_k\ra {\bf G}_k$ is an open map.
\item There exists a basis of compact neighborhoods $W$ of $e$ in ${\bf P}_k$
such that
$$a^{-1}Wa\subset W \quad\text {for all $a\in {\bf A}_k^+$}.$$
\item For every neighborhood $U$ of $e$ in $G$, there exists
a neighborhood $V$ of $e$ in $G$
such that $$HVa\subset HaU   \quad\text {for all $a\in {\bf A}_k^+$}.$$
\end{enumerate}
\el

\proof
(1) When $\text{char}(k)=0$, it follows
from the fact that
Lie algebras of ${\bf P}_k$ and ${\bf H}_k$ generate the Lie algebra
of ${\bf G}_k$ as a vector space.
For a characteristic free argument,
see \cite{HW}
or Proposition I.2.5.4 in \cite{Ma2}.

(2) When $\text{char}(k)=0$,  note that
the action of ${\bf A}_k^+$
on the Lie algebra
$\text{Lie}({\bf P}_k)$ gives a family of commuting
semisimple linear maps ${\rm Ad}\,(a)$ whose eigenvalues
have bounded above by $1$ in their absolute values.
It follows that there exists a basis of compact neighborhoods $W_0$ of
$0$ in  $\text{Lie}({\bf P}_k)$
which are invariant by all ${\rm Ad}\,(a)$, $a\in {\bf A}_k^+$.
 It suffices to set $W={\rm exp}(W_0)$.

It is easy to adapt this argument
in positive characteristic case; write
${\bf P}_k={\bf L}_k{\bf N}_k$, and
note that ${\bf L}_k$ contains a
${\bf A}_k$-invariant compact open subgroup.
Now considering the linear group
${\bf N}_k$ as a group of
upper triangular matrices in a suitable
basis where elements of ${\bf A}_k^+$ are diagonals
with increasing coefficients in absolute value,
we can find a basis of compact neighborhoods
$W$ as desired.

(3) Choose $W$ as in (2) small enough so that $W\subset U$ and
choose any neighborhood
$V$ of $e$ in $G$ contained in $HW$.
We then have $HVa\subset HWa\subset HaU$, as required.
\hb\vs

\noindent {\it Proof of Proposition \ref{prowavefront}. }
We will prove that $G$ has the wavefront property on $H\ba G$
with the subset $F=AK$ defined in Lemma \ref{lemkah}.
By Lemma \ref{lemkah}, it suffices to show that
{\it for every neighborhood $U$ of $e$ in $G$, there exists
a neighborhood $V$ of $e$ in $G$ such that
 $HVg\subset HgU$ for all $g\in AK$}

Recall $A=\cup_{i=1}^m {\bf A}_i(k)$ where
${\bf A}_i$ is a maximal $(k,\si)$-split torus of $G$.
Fix $i$ and a positive Weyl chamber $C$ of
${\bf A}_i(k)$.

Since $K$ is a compact set, there exists a neighborhood $U_0$
of $e$ in $G$ such that $k^{-1}U_0k$ is contained in $U$
for all $k$ in $K$.
By Lemma \ref{lemawaw}.(3),
there exists a neighborhood $V_C$ of $e$ in $G$ such that
 $V_Ca\subset HaU_0 $ for all $a\in C$.

Now set $V:=\cap_{C}V_C$ where the intersection is taken
over all (finitely many) positive Weyl chambers of
${\bf A}_i(k)$, $1\le i\le m$.
Then for $g=ak\in (\cup_C C)K=AK$ with $k\in K$ and $a\in C$, we have

$$
HVg \subset
 HV_C ak \subset
HaU_0k  \subset
HakU = HgU\; .
$$
This finishes the proof.

In section \ref{secapplication}, we will use
this wavefront property in the product situation,
owing to the following straightforward proposition.

\bp
\label{prowavefront2}
For $i=1,\ldots ,m$, let $G_i$ be a locally compact group,
$H_i\subset G_i$ a closed subgroup,
$G:=\prod_{i=1}^m G_i$ and $H:=\prod_{i=1}^m
H_i$.
If $G_i$ has the wavefront property on $H_i\ba G_i$
for each $1\le i\le m$,
then $G$ has the wavefront property on $H\ba G$.
\ep

The following theorem is an immediate consequence of
Propositions \ref{promixing}, \ref{prohowemoore},
\ref{prowavefront} and \ref{prowavefront2}.

\bt
\label{thwavefront}
For $i=1,..,m$, let $k_i$ be a local field, $G_i$ the group
of $k_i$-points of a semisimple $k_i$-group, $\si_i$ an
involution of $G_i$ defined over $k_i$,
$G_i^{\si_i}$ its group of
fixed points
and $H_i$ a closed subgroup of finite index of $G_i^{\si_i}$.
Let $G=\prod_{i=1}^m G_i$ and $H:=\prod_{i=1}^m
H_i$.

Then the group $G$ has the wavefront property on $H\ba G$.

Moreover, for any irreducible lattice $\Ga$ in $G$,
the action of $G$ on $\Ga\ba G$ is mixing.
\et

Note that this theorem provides
many natural examples of
 triples $(G, H, \Gamma)$
which satisfy the hypothesis of the propositions
\ref{proequihorbit}, \ref{prosumintegral} and
\ref{prowellrounded}.

\section{Equidistribution of translates of $H$-orbits}
\label{sechorbits}
In this section,
let $G$ be a locally compact group, $H \subset G$ a closed subgroup,
$\Ga\subset G$ a lattice such that $\Ga_H:=\Ga\cap H$ is a lattice in
$H$. Set $X=\Gamma\ba G$ and $Y=\Ga_H\ba H$.
Let $\mu_X$ and $\mu_Y$ be invariant measures on $X$ and $Y$
respectively.

\bp
\label{proequihorbit}
Suppose that the action of $G$ on $X$ is mixing and
that $G$ has the wavefront property on $H\ba G$.
Then the translates $Yg$
become equidistributed
in $X$, as $g\to \infty$ in $H\ba G$.
\ep

This means that as the image of $g$ in $H\ba G$ leaves every
compact subsets, the sequence of probability measures
$\frac{1}{\mu_Y(Y)}g_*\mu_Y$ weakly
converges to $\frac{1}{\mu_X(X)}\mu_X $,
i.e., for any $ \psi\in C_c(X)$, we have
\begin{eqnarray}
\label{eqnlaylax}
\frac{1}{\mu_Y(Y)}\int_Y \psi (yg) d\mu_Y(y)\to
\frac{1}{\mu_X(X)}\int_X \psi\; d\mu_X \; .
\end{eqnarray}

\proof
The following proof is adapted from \cite{EM}; we point out that
the case when $Y$ is non-compact requires a bit more care, which
was not addressed in \cite{EM}.
Since $G=HF$, we may assume that $g$ belongs to the subset $F$ in
Definition \ref{defwavefront}.
We assume, without loss of generalities, that $\mu_X$ and
$\mu_Y$ are probability measures.
Let $\psi\in C_c(X)$.
Fix $\eps>0$.
By the uniform continuity of $\psi$ there exists a
neighborhood $U$ of $e$ in $G$ such that
\begin{eqnarray}
\label{eqnequic}
|\psi(xu)-\psi(x)|<\eps\; \quad\text{ for all $u\in U$ and
$x\in X$}.
\end{eqnarray}
By the wavefront property of $G$ on $H\ba G$,
there exists a compact neighborhood $V\subset U$ of $e$
in $G$ such that
\begin{eqnarray}
\label{eqngvugh}
Vg\subset HgU \quad \text{for all } g\in F
\end{eqnarray}

Choose a compact subset $Y_\epsilon \subset Y$ of measure at
least $\mu_Y(Y_\epsilon)\geq 1-\eps$.
Choose a Borel subset $W\subset V$ in $G$ transversal  to $H$,
i.e., a subset $W$ of $G$ such that the multiplication
$m:H\times W \ra G$ is injective with the image
$HW$ being an open neighborhood of $e$ in $G$.
 Using the compactness of $Y_\epsilon$ and the discreteness of
$\Gamma$,
we may assume that the image of $W$ in $H\ba G$ is small enough so that the
multiplication
$m:  Y_\epsilon\times W \ra Y_\epsilon W$ is
a bijection\footnote{When $Y$ is compact, one can choose
the transversal $W$
such that the map $Y\times W\ra YW$ is bijective
onto an open subset of $X$.
When $Y$ is not compact, such a transversal does not always exist.
Here is an example: let $G$ be the orthogonal group
of the quadratic form $x^2+y^2+z^2-t^2$ on $\m R^4$,
$v_0=(1,0,0,0)$, $v_1=(1,0,2,2)$,
$\Ga= G_{\m Z}$ and $H$ the stabilizer of the point $v_0$.
One checks easily that (a) $v_1=\ga v_0$ for some $\ga\in \Ga$
and  that (b) $v_0$ is a limit of elements $v_1h_n$
of the $H$-orbit of $v_1$.
Hence there exists a sequence $g_n$ converging to $e$ in $G$
such that $H g_n\cap \gamma H\not = \emptyset$.\\
To check (a), take $\ga=
\left(\begin{array}{cccc}
1 & 0 & 2 & 2 \\
0 & 1 & 0 & 0 \\
2 & 0 & 1 & 2 \\
2 & 0 & 2 & 3
\end{array}
\right)$.
For (b), take $h_n=
\left(\begin{array}{cccc}
1 & 0 & 0 & 0 \\
0 & 1 & 0 & 0 \\
0 & 0 & \mbox{ }\cosh n & -\sinh n \\
0 & 0 & -\sinh n & \mbox{ }\cosh n
\end{array}
\right)$.
}
onto its image $Y_\e W\subset X$.

Let $\mu_W$ be the measure on $W$ such that
$\mu_X=\mu_Y\mu_W$ locally.

Setting
$$I_g:=\int_Y\psi (yg)d\mu_Y(y),$$
we need to show that
\begin{eqnarray}
\label{eqnig}
I_g \to \int_{X}\psi\; d\mu_X
\quad\text{as $g \in F$ goes to infinity in $G$.}\end{eqnarray}

For simplicity,
set
$$
J_g={\textstyle\frac{1}{\mu_W(W)}}\int_{Y\times W}\!\!\!
\psi (ywg)\, d\mu_Y(y)\,
d\mu_W(w)
\;\; {\rm and}\;\;
K_g={\textstyle\frac{1}{\mu_W(W)}}\int_{Y_\e W}
\!\!\psi (xg)\, d\mu_X(x).$$

Roughly speaking, we will argue that $I_g$
is close to $J_g$ as a consequence
of the wavefront property,  $J_g$ is close to $K_g$
since the volume of $Y\! -\! Y_\e$ is small, and
finally $K_g$ is close to the average of $\psi$ for large $g$
 because of the
mixing property.

By (\ref{eqngvugh}),
for each $w\in W$ and $g\in F$,
we have $wg=h_{g, w}gu$ for some
$h_{g, w}\in H$. Hence
\begin{eqnarray*}
\label{eqnig2}
|I_g -
\int_Y \psi (ywg)\, d\mu_Y(y)|
&=& |\int_Y \psi (yg)\, d\mu_Y(y) -
\int_Y \psi (y h_{g, w}gu)\, d\mu_Y(y)|
\\
&=&|\int_Y (\psi (yg)-\psi (ygu))\, d\mu_Y(y)|
\le \e \quad\text{by (\ref{eqnequic})}.
\end{eqnarray*}

Therefore we have
$$|I_g-J_g|\le \e .$$

By the choice of $W$, we have
$$K_g=\frac{1}{\mu_W(W)} \int_W\int_{Y_\e} \psi (ywg)\,
d\mu_Y(y)\, d\mu_W(w)$$
and hence
$$|J_g-K_g| \le 2\mu_Y(Y-Y_\e)\,  \|\psi\|_\infty
\le 2 \|\psi\|_\infty \,\eps.$$
 Since $K_g=
\frac{1}{\mu_W(W)}\int_{X}\psi(xg)\; {\bf 1}_{Y_\e W} (x)d\mu_X(x)$
where ${\bf 1}_{Y_\e W}$ is the characteristic function of $WY_\e$,
the mixing property of $G$ on $\G\ba G$ says that
$K_g$ converges to $\mu_Y (Y_\e)\int_X\psi d\mu_X $ as $g \to\infty$ in $F$.
Hence for $g\in F$ large enough
we have,

\begin{align*}
|K_g-\int_X \psi \, d\mu_X|
&\le \eps +\mu_Y(Y-Y_\eps)\int_X\psi \,d\mu_X
\le (1+\| \psi\|_\infty )\,\eps.
\end{align*}

Putting this together, we get
\begin{align*}
|I_g-\int_X \psi \, d\mu_X|& \le
|I_g-J_g| + |J_g-K_g| +
|K_g- \int_X\psi \,d\mu_X |\\
&\le
 (2+3\|\psi\|_\infty)\, \e .
\end{align*}

Since $\e>0$ is arbitrary,
this shows the claim.
\hb\vs

Using Theorem \ref{thwavefront}, we obtain:

\bc
\label{corequihorbit}
Let $G$, $H$, $\Ga$ be as in Theorem \ref{thwavefront}.
Then the translates $Yg:=\Ga_H\ba H g$
become equidistributed
in $X:=\Ga\ba G$ as $g\to\infty$ in $H\ba G$.
\ec

\section{Sums and integrals}
\label{secsumintegral}
Let $G$ be a locally compact group, $H\subset G$ a closed subgroup,
$\Ga\subset G$ a lattice such that $\Ga_H:=\Ga\cap H$ is a lattice in
$H$. Let $x_0:=[\Ga]$ be the base point in $X:=\Ga\ba G$, $Y=x_0H$
and $z_0:=[H]$ be the base point in $Z:=H\ba G$. We note that
$z_0\Gamma$ is a discrete subset of $Z$.
There exist  $G$-invariant measures $\mu_X$, $\mu_Y$ and $\mu_Z$
on $X$, $Y$ and $Z$.
We normalize them so that $\mu_X=\mu_Z\mu_Y$ locally.

For a given sequence  of non-negative functions $\ph_n$ on $Z$
with compact support,
we define  a function $F_n$ on $X$ so that,  for $x=x_0 g$,
$F_n(x)$ is the sum of $\ph_n$ over the discrete
orbit $z_0\Gamma g$:

\begin{eqnarray}
\label{eqnfn}
F_n(x):=\sum_{\ga\in \Ga_H\ba \Ga}\ph_n(z_0\ga g)\quad\text{for $x=\Ga g$}.
\end{eqnarray}
We would like to
 compare the values of $F_n$ with the space average
over $Z$:
\begin{eqnarray}
\label{eqnin}
I_n:=\tfrac{\mu_Y(Y)}{\mu_X(X)}\int_Z\ph_n(z) \;d\mu_Z(z)
\end{eqnarray}
We remark that this normalized integral $I_n$ does not depend
on the choices of measures.

The following proposition \ref{prosumintegral}
says that the sum $F_n$  is asymptotic
to the normalized integral $I_n$, at least weakly.

\bp
\label{prosumintegral}
Suppose that the translates $Yg$
become equidistributed
in $X$ as $g\to\infty$ in $Z$.
Then for any sequence  of non-negative functions $\ph_n$ on $Z$
with compact support such that $\max_n \|\ph_n \|_\infty <\infty$ and
${\D\lim_{n\ra \infty}}\int_Z\ph_n \, d\mu_Z =\infty$,
the ratios $F_n(x)/I_n$ converge weakly to $1$
as $n \to\infty$.
\ep

This means that, for all $\al\in C_c(X)\; ,\;\;\;$
\begin{eqnarray}
\label{eqnsumintegral}
\lim_{n\ra\infty}\frac{1}{I_n}\int_X F_n(x)\;\al(x)d\mu_X(x)
&=&
\int_X\al (x) d\mu_X(x)\; .
\end{eqnarray}
\vs

\proof
Using transitivity properties
for invariant integration on homogeneous spaces,
we obtain that
for all $\al\in C_c(X)$,
\begin{eqnarray*}
\int_{\Gamma\ba G} F_n\al d\mu_{\Ga \ba G}
&=&
\int_{\Ga \ba G} \sum_{\ga\in \Ga_H\ba \Ga}\ph_n( H \ga g)\al(\Ga g)
d\mu_{\Ga \ba G}(\Ga g)\\
&=&
\int_{ G_H\ba \Ga} \ph_n(Hg)\al(\Ga g)
d\mu_{\Ga \ba G_H}(\Ga_H g)\\
&=&
\int_{H\ba G} \int_{\Ga_H\ba H}\ph_n(Hg)\al(\Ga_H h g)
d\mu_{\Ga_H\ba H}(\Ga_H h)d\mu_{H\ba G}(H g)\\
&=&
\int_{H\ba G} \ph_n(z)\be(z)d\mu_{H\ba G}(z)
\end{eqnarray*}
where $\be$ is the function on $Z$ given by,
\begin{eqnarray*}
\be(Hg)
&=&\int_{\Ga_H \ba H}\al(\Ga_H h g )d\mu_{\Ga_H\ba H}(\Ga_H h)\\
&=&
\int_Y\al(yg)d\mu_Y(y)\; .
\end{eqnarray*}
By assumption, we have
\begin{eqnarray*}
\D{}\;\;\;
\lim_{z\ra\infty} \be(z) =
\tfrac{\mu_Y(Y)}{\mu_X(X)}\int_X \al (x) d\mu_X(x)\; .
\end{eqnarray*}
Since $I_n=\int_Z\ph_n \ra \infty$ and $\ph_n$ are uniformly bounded,
by the dominated convergence theorem
\begin{eqnarray*}
\lim_{n\ra\infty}\frac{1}{I_n}\int_{Z} \ph_n(z)\be(z)d\mu_{Z}(z)
= \int_X \al \;d\mu_X .
\end{eqnarray*}
Hence we obtain the equality \eqref{eqnsumintegral}.
\hb

\section{Counting and equidistribution }
\label{seccounting}
We will now improve the weak convergence in proposition
\ref{prosumintegral} to the
 pointwise convergence of the functions $F_n$.
This requires some hypothesis on the sequence  of
functions $\ph_n$
which will be called {\it well-roundedness}.
We keep the notations of section \ref{secsumintegral}.

\bd
\label{defwellrounded}
A sequence or a family of non-negative integrable functions
$\ph_n$ of $Z$ with compact support
 is said to be {\it well-rounded} if for any $\eps>0$,
there exists a neighborhood $U$ of $e$ in $G$,
such that the following holds for all $n$.
\begin{eqnarray*}
\label{eqnwellrounded}
(1\! -\! \eps)\int_Z \; (\sup_{u\in U}\ph_n(zu)) \; d\mu_Z(z)\leq
\int_Z \ph_n\; d\mu_Z
\leq
(1\! +\! \eps)\int_Z \; (\inf_{u\in U}\ph_n(zu))\; d\mu_Z(z).
\end{eqnarray*}

A sequence  of subsets $B_n$ of $Z$ is said to be well-rounded if
the sequence ${\bf 1}_{B_n}$ is well-rounded.
\ed

Sometimes we will apply the above definition to a
continuous family $\{\ph_T\}$ of functions or subsets, whose meaning
should be clear.

Recall that we want to compare
the orbital sum
$F_n(x_0)=\sum_{\ga\in \Ga/\Ga_H}\ph_n(\ga z_0)$ with
the average
$I_n=\tfrac{\mu_Y(Y)}{\mu_X(X)}\int_Z\ph_n(z) d\mu_Z(z)$.

\bp
\label{prowellrounded}
Keep the notations and hypothesis
of Proposition \ref{prosumintegral}, and
assume that the sequence $\ph_n$ is well-rounded.
Then,
$$F_n(x_0)\sim I_n\quad\text{as $n\to \infty$} .$$
\ep

The notation $a_n \sim b_n$ means that the ratio of $a_n$ and $b_n$
tends to $1$ as $n\to \infty$.
\vs

\proof
Once again, we may normalize the measures so that
$\mu_X(X)=\mu_Y(Y)=1$.
Fix $\eps>0$ and let $U$ be a neighborhood
of $e$ in $G$
given by Definition \ref{defwellrounded}.
We introduce the functions $\ph_n^\pm$ on $Z$ defined by
\begin{eqnarray*}
\label{eqnwellrounded2}
\ph_n^+(z):=\sup_{u\in U}\ph_n(zu^{-1})
& {\rm and}&
\ph_n^-(z):=\inf_{u\in U}\ph_n(zu^{-1})
\end{eqnarray*}
and their integrals $I_n^\pm:=\int_Z\ph_n^\pm \; d\mu_Z$.
Note that for each $n$,
\begin{eqnarray}
\label{eqnininin}
(1-\eps)I_n^+
\leq I_n
\leq
(1+\eps)I_n^-\; .
\end{eqnarray}
We also introduce the functions $F_n^{\pm}$ on $X$:
\begin{eqnarray*}
\label{eqnfnpm}
F_n^\pm(x)=\sum_{\ga\in \Ga_H\ba \Ga}\ph_n^\pm( z_0 \ga g)\;
\quad\text{for } x=\Ga g .
\end{eqnarray*}

It is easy to check that,
for all $u\in U$ and $x\in X$
\begin{eqnarray*}
\label{eqnfnfnfn}
F_n^-(ux)
\leq
F_n(x)
\leq
F_n^+(ux)
\; .
\end{eqnarray*}

Choose a non-negative continuous function $\al$
on $X$ with $\int_X\al =1$ and
with support included in $x_0 U$
so that the following holds for all $n$:
\begin{eqnarray*}
\label{eqnfnal}
\int_X \al\, F_n^- \, d\mu_X
\leq
F_n(x_0)
\leq
\int_X \al \, F_n^+\, d\mu_X
\; .
\end{eqnarray*}
Applying Proposition \ref{prosumintegral}
to the sequences of functions $\ph_n^\pm$,
we obtain, for all $n$ large,
\begin{eqnarray}
\label{eqninfnin}
(1-\eps)I_n^-
\leq
F_n(x_0)
\leq
(1+\eps)I_n^+
\; .
\end{eqnarray}
Using the estimations (\ref{eqnininin}) and (\ref{eqninfnin}),
every cluster value of the sequence of ratios
$F_n(x_0)/I_n$ is within the interval
$[\frac{1-\eps}{1+\eps},\frac{1+\eps}{1-\eps}]$.
Hence this sequence converges to $1$.
\hb

\section{Well-roundedness}
\label{secwellrounded}
In this section, we provide explicit examples of
well-rounded sequences $\ph_n$
in order to apply Proposition \ref{prowellrounded}.
We start with an observation that the product of
 well-rounded sequences is again well-rounded.

\begin{Ex}
\label{ex1}
For each $i=1,\ldots ,m$, let $G_i$ be a locally compact
 group,
$H_i\subset G_i$ a closed subgroup
and $\ph_{i,n}$ be a well-rounded sequence of functions on
$Z_i:=H_i\ba G_i$.
Let $G:=\prod_{i=1}^m G_i$, $H:=\prod_{i=1}^m H_i$ and $Z:=\prod_{i=1}^m
 Z_i$.
Then the sequence $\ph_n$ defined by
$$\ph_n(z_1,\ldots ,z_m)=\prod_{1\leq i\leq m} \ph_{i,n}(z_i)$$
is well-rounded.
\end{Ex}

\proof
Fix $\eps>0$. Let $U_i$ be a neighborhood of $e$ in $G_i$
such that the functions $\ph_{i,n}^\pm$ on $Z_i$
defined by
$\D\ph_{i,n}^+(z_i)=\sup_{u_i\in U_i}\ph_{i,n}(z_i u_i)$
and
$\D\ph_{i,n}^-(z_i)=\inf_{u_i\in U_i}\ph_{i,n}(z_iu_i)$
satisfy
\begin{eqnarray*}
\label{eqnphnipm}
\textstyle
(1-\eps)\int_{Z_i}\ph_{i,n}^+
\leq
\int_{Z_i}\ph_{i,n}
\leq
(1+\eps)\int_{Z_i}\ph_{i,n}^-
\; .
\end{eqnarray*}
Let $U:=\prod U_i$, and $\ph_{n}^\pm$
be the functions on $Z$
defined by
$$\D\ph_{n}^+(z)=\sup_{u\in U}\ph_{n}(zu)\;\quad\text{and}\quad
\D\ph_{n}^-(z)=\inf_{u\in U}\ph_{n}(zu)$$
so that $\ph_n^\pm =\prod \ph_{i,n}^\pm$ and
\begin{eqnarray*}
\label{eqnphnpm}
\textstyle
(1-\eps)^m\int_{Z}\ph_{n}^+
\leq
\int_{Z}\ph_{n}
\leq
(1+\eps)^m\int_{Z}\ph_{n}^-
\; .
\end{eqnarray*}
Hence the sequence $\ph_n$ is well-rounded.
\hb\vs

The next example deals with the constant sequences.
It will be used both for the archimedean
and the non-archimedean factors.

\begin{Ex}
\label{ex2}
Let $G$ be a locally compact group, $H$ a closed
 subgroup of
$G$, $Z=H\ba G$, and $\ph\in C_c(Z)$ with $\ph \ge 0$ and $\ph
\neq 0$.
Then the constant sequence
$\ph_n=\ph$ is well-rounded.
\end{Ex}

\proof  Use the uniform continuity of $\ph$ and
the compactness of its support.
\hb\vs

The following example
 deals with the archimedean factors.
\begin{Ex}
\label{ex3}
Let $G$ be a  real  semisimple Lie group
with finitely many connected components,
$V$ a finite dimensional representation of $G$,
$Z$ a closed $G$-orbit in $V$ with
an invariant measure $\mu$
and $\| .\|$ an euclidean norm on $V$.
Then the family of balls $B_T:=\{ z\in Z\mid \| z\|\leq T\}$, $T\gg 1$
is well-rounded.
\end{Ex}

\proof  By Corollary \ref{cordevas}.a of the appendix,
we have $$\mu(B_T)\sim_T c\, T^a(\log T)^{b}$$
for some  $a\in \q_{\ge 0}$,  $b\in \z_{\ge 0}$ and $c>0$.
It is easy to deduce the claim from the above asymptotic
using the assumption that the action of $G$ is linear on $Z$.
\hb\vs

As for the non-archimedean factors, we have:
\begin{Ex}
\label{ex4}
Let $k$ be a non-archimedean local field,
$G$ the group of $k$-points
of a connected semisimple $k$-group,
$\rho:G\ra GL(V)$ a representation of $G$ defined over $k$,
$Z$ a closed $G$-orbit in $V$ with
an invariant measure  and
$\|\cdot \|$ a norm on $V$.
Then, both the family
of non empty balls $B_T:=\{ z\in Z\;|\; \| z\| \leq T\}$,
and the family of
non-empty spheres $S_{T}:=\{ z\in Z\;|\;\| z\| =T\}$,
 are well-rounded.
\end{Ex}

\proof
Since the action of $G$ on $V$ is linear,
the stabilizer in $G$ of the norm  is a compact open subgroup of $G$.
Hence this example is a special case of the
following easy
assertion.
\hb\vs

\begin{Ex}
\label{easy}
\label{ex5}
Let $G$ be a locally compact (unimodular) group,
$H$ a closed (unimodular) subgroup of
$G$, $Z=H\ba G$  and $U$ a compact open subgroup
of $G$. Then any sequence $\ph_n$ of non-negative $U$-invariant
$L^1$-functions on $Z$ is well-rounded.\end{Ex}

As the last example, we will show that a sequence
of the height balls is well rounded.
We will need the following basic lemma.
\bl
\label{lemarithm}
Let $I$ be a finite set. For each $i\in I$, let $\tau_i> 1$
and $\la_i>0$ be given.
Let $\la:\m N^I\to \m R^+$ be given by
$\la(m)= \sum_{i\in I}\la_i m_i $ for $m=(m_i)$, and
 $P:\m N^I \to \m R^+$ a function given by
$P(m)=\prod_{i\in I} P_i(m_i)$  where
$P_i$ is a real-valued function of a variable $x$ given by a polynomial
expression in $(x,\tau_i^x)$ and which is positive on $\m N$.
Then we have
 $$w_{t+1}=O(w_t) \quad\text{for $t$ large}$$
where
$$
w_t:=\sum_{\{ m\in \m N^I,\,\la(m)\leq t\} }P(m)\, .$$
\el

\proof
Since each $P_i$ is positive on $\m N$,
there exists $C>0$ such that for $i\in I$ and $n\in \m N$,
one has
$$
P_i(n+1)\leq C\, P_i(n)\, .
$$
Hence for  each $m\in \m N^I$
and each $e$ in the basis $E$ of $\m N^I$,
one has
$$
P(m+e)\leq C\, P(m)\, .
$$
Setting $t_0:=\min_i \la_i = \min_{e\in E} \lambda(e)$,
one gets
$$
w_{t+t_0}\leq
\sum_{e\in E}\sum_{\la(m)\leq t}P(m+e)
\leq rCw_t\, .
$$
for $r=|I|$.
Hence we conclude that $w_{t+1}\leq (rC)^kw_t$ with $k =\frac{1}{t_0}+1$.
\hb
\vs

\br
\label{remwt}
{\rm
One can improve the conclusion of Lemma \ref{lemarithm}:}
there exist $a\geq 0$, $b\in \m Z_{\geq 0}$
and $c_1,c_2>0$ such that, for $t$ large,
$$
c_1 e^{at}t^b
\leq w_t \leq
c_2 e^{at}t^b\; .
$$
Moreover, setting $C_i\tau_i^{d_ix}x^{b_i}$
for the dominant term of $P_i(x)$,
the exponents $a$ and $b$ are given 
respectively by $e^a=\D\max_{i\in I}\tau_i^{d_i/\la_i}$
and $b$ is given by $b=\sum_i (b_i +1)-1$ where the sum is 
taken over all $i$ such that $e^a=\tau_i^{d_i/\la_i}$.
{\rm The proof is a straightforward induction on $|I|$. 
Here is a sketch:
one may assume that, for all $i$, $P_i(x)=C_i e^{a_ix}x^{b_i}$
and $\la_i=1$. One fixes 
$i_0\in I$, set $\check{I}:=I-\{ i_0\}$ 
and writes
$w_t:=\sum_{1\leq n\leq t}P_{i_0}(n)\check{w}_{t-n}$
where, by induction hypothesis, 
$\check{w}_t$ satisfies a similar estimation, as $t\to\infty$:
$
\check{c}_1 e^{\check{a}t}t^{\check{b}}
\leq \check{w}_t \leq
\check{c}_2 e^{\check{a}t}t^{\check{b}}\; 
$
for some $\check{a}\geq 0$, $\check{b}\in \m Z_{\geq 0}$
and $\check{c}_1,\check{c}_2>0$.
From that, one gets the required estimation for $w_t$.
}
\er

For the rest of this section,
let $I$ be a finite set.
For each $i\in I$, let $k_i$ be a local field
of characteristic $0$,
$G_i$  the group
of $k_i$-points of an algebraic $k_i$-group,
$V_i$ an algebraic representation of $G_i$,
 and $Z_i\subset V_i$ a non-zero closed $G_i$-orbit
with an invariant measure $\mu_i$.

We set
$$G:=\prod_{i\in I}G_i, \;Z:=\prod_{i\in I} Z_i,\;\mu:=
\otimes_{i\in I}\mu_i.$$

Let $I_\infty\subset I$ be the set of
indices with $k_i$ archimedean, and $I_f:=I\setminus I_\infty$.
The partition $I=I_\infty \sqcup I_f$ induces decompositions
$G=G_\infty\times G_f$ of the group,
$Z=Z_\infty\times Z_f$ of the orbit, and
$\mu=\mu_\infty\otimes \mu_f$ of the invariant measure.

Let $\|\cdot \|_i$ be a norm on $V_i$. We assume that
$\|\cdot\|_i$ is euclidean if $i\in I_\infty$ and a max norm otherwise.
These norms define a height function $h: Z \to \m R^+$
$$h(z)=\prod_{i\in I} \|z_i\|_i .$$
Since each $Z_i$ is a closed non-zero subset in $V_i$,
we have $\min_{z\in Z_i} \|z\|_i >0$ and hence
$h$ is a proper function on $Z$. We can also write $h=h_\infty \otimes h_f$
where $h_\infty:=\prod_{i\in I_\infty} \|\cdot \|_i$
and $h_f:=\prod_{i\in I_f} \|\cdot \|_i$.

Set
$$
b_T:=\{z\in Z_\infty \mid  h_\infty (z)\leq T\}
\;\; {\rm ,}\;\;
\be_T:=\{z\in Z_f \mid  h_f (z)\leq T\}\, ,$$
$$
V_T:=\mu(B_T)
\;\; ,\;\;
v_T=\mu_\infty(b_T)
\;\; {\rm and}\;\;
w_T=\mu_f(\be_T)\; .
$$

\bl
\label{lemex6}
Assume that $h_\infty$ is not constant on $Z_\infty$.
\begin{enumerate}
\item There exist  $a\in \m Q_{\ge 0}$, $b\in \m Z_{\ge 0}$,
and $c>0$ such that
 as $T\to \infty$, $$
v_T
\sim c \, T^a (\log T)^{b}
\;\;\;\mbox{ and}\;\;\;
\tfrac{d}{dT}v_T
\sim c \tfrac{d}{dT}( T^{a} (\log T)^{b})
.$$
\item There exist constants $\kappa>0$ and $C_1>0$ such that
for all $\eps\in ]0,1[$ and all $T\geq 0$
$$v_{(1+\eps)T}- v_T\leq C_1(v_T+1)\eps^\kappa .$$
\item For $T$ large, one has $w_{2T}= O(w_T)$.
\item There exist $\kappa>0$ such that
uniformly for $T$ large and $\eps\in ]0,1[$,
$$V_{(1+\eps)T}=(1+O(\eps^\kappa))V_T\, .$$\end{enumerate}
\el

\proof
(1): Apply Proposition \ref{prodevas} of the appendix
to the regular function $F:=h_\infty^2$
on the orbit $Z_\infty$. Note that since $v_T$ is an increasing function
of $T$, one has $a\geq 0$.
Moreover, note that, when $a=b=0$, the orbit is of finite volume
and hence compact.

(2): First note that, since $h_\infty$ is not constant on $Z_\infty$,
the function $v_T$ is continuous.

When $T$ is large, we use (1) to get
the following bound
$v_{(1+\eps)T}-v_T =O(\eps v_T)$
which is uniform in $\eps\in ]0,1[$.

When $T$ is bounded, we use the fact that the function
$v_T$ is differentiable
except at the critical values $\tau$ of $h_\infty$.
Since $h_\infty^2$ is a regular function,
there are only finitely many such critical values $\tau$.
Around these points,
there exists a constant $\kappa$, $0<\kappa <1$ such that,
for $\eps>0$ small,
one has the following bound for the derivative:
\begin{eqnarray}
\label{eqnvptau}
v'_{\tau \pm\eps}
&=&
O(\eps^{\kappa-1}).
\end{eqnarray}
This assertion is a consequence of
Theorem \ref{thdevas} of the appendix.
More precisely, set $f:=\pm(h_\infty -\tau)$.
Since $Z_\infty$ is smooth, one can choose
$\eps_0>0$ and an open covering $U_j$ of $f^{-1}(]-\eps_0,\eps_0[)$
by open sets bianalytically homeomorphic to balls.
A partition of unity gives us
C$^\infty$ functions $\ph_j$ with compact support in $U_j$
such that $\sum_j\ph_j=1$ on $f^{-1}(]-\eps_0/2,\eps_0/2[)$.
We simply apply Theorem \ref{thdevas} to these functions
$f$ and $\ph_j$ to get \eqref{eqnvptau}.

Integrating $v'_t$ on the interval
$[T,(1+\eps)T]$, and using \eqref{eqnvptau}
near the critical values in this interval,
one gets, uniformly for $\eps$ small and $T$ bounded,
$$v_{(1+\eps)T}-v_T =O(\eps^\kappa).$$

Putting these together proves the claim.

(3):
We will assume, as we may, that $\inf_{z\in Z_i}\|z\|_i\geq 1$
for each $i\in I_f$.
For any $|I_f|$-tuple
$m=(m_i)_{i\in I_f}\in \m N^{I_f}$,
we set
$$S(m)=\prod_{i\in I_f}S_i(m_i)\quad\text{where}
\quad S_i(m_i):=\{ z\in Z_i : \| z\|_i =m_i\} .$$

Letting
$\om_m:=\mu_f(S(m))$ and
$\pi_m:=\prod_i m_i$, one has
$$
\textstyle w_T =\sum_{\{m\in \c M,\pi_m\leq T\}}\om_m
$$
where $\c M\subset \m N^{I_f}$
consists of $m\in \m N^{I_f}$ with non-empty $S(m)$.
 The main point of the proof is to use the formula for
$\om_m$ given by Theorem \ref{thden} of the appendix.
According to this formula, there is a finite partition of $\c M$
in finitely many pieces $\c M_\al$ such that\\
- each piece $\c M_\al$ is a product of subsets $\c M_{\al,i}$ of $\m N$
which are either points or of the form
$\{ m_i=c_{\al,i}\, q_{\al,i}^{n_i}: n_i\in \m N\}$
for some positive integers $c_{\al,i}, q_{\al,i}$,\\
- on each piece $\c M_\al$, the volume $\om_m$ is given by
a formula $\prod_{i\in I_f} P_{\al,i}(n_i,q_{\al,i}^{n_i/d})$
where $P_{\al,i}$
is a polynomial and $d$ a positive integer.

According to Lemma \ref{lemarithm} with $T=2^t$,
the volume
$$
\textstyle w_{\al,T} :=\sum_{\{m\in \c M_\al,\pi_m\leq T\}}\om_m
$$
satisfy the bound $w_{\al,2T}=O(w_{\al,T})$.
Hence one has $w_{2T}=O(w_T)$ as required.

(4): Let $T_0:=\inf_{z\in Z_\infty}h_\infty(z)>0$.
According to (2), there exists  $C>0$
such that for $T$ large
\begin{eqnarray*}
V_{(1+\eps)T} -V_T&=&
\textstyle\sum_{m} (v_{(1+\eps)T/\pi_m}-v_{T/\pi_m})\, \om_m\\
&\leq &
\textstyle
C \eps^\kappa\left(
\sum_{m} v_{T/\pi_m}\, \om_m
+\sum_{m} \om_m \right)\\
&\leq &
C \eps^\kappa(V_T+w_{2T/T_0})
\end{eqnarray*}
where the above sums are over all the multi-indices $m\in \m N^{I_f}$
with $T_0\pi_m\leq 2T$.
Then, applying (3) twice,  there exists $C'>0$
such that for $T$ large
\begin{eqnarray*}
V_{(1+\eps)T} -V_T &\leq &
C \eps^\kappa(V_T+C'w_{T/2T_0})\\
&\leq &
C(1+C'v_{2T_0}^{-1})\eps^\kappa V_T,
\end{eqnarray*}
as required, since $v_{2T_0}>0$.
\hb

\br
\label{remvt} 
{\rm 
One has the following estimate for the
volume $V_T$ of the height ball :}
there exist $a\in \m Q_{\geq 0}$, $b\in \m Z_{\geq 0}$
and $c_1,c_2>0$ such that for $T$ large,
$$
c_1 T^a\log (T)^b
\leq V_T \leq
c_2 T^a \log (T)^b\; .
$$
{\rm 
This is a straightforward consequence of 
the formula
$V_T:=\int_0^\infty w_{T/t} v'_{t}dt$ and of the estimation of $w_T$ and  $v'_T$
given in Remark \ref{remwt} and Lemma \ref{lemex6} (1). 
}
\er

\bp
\label{prohbwell}
\label{proex6}
{\bf (Height ball)}
The family of height balls
$B_T:=\{z\in Z \mid  h (z)\leq T\}$, $T\gg 1$,
is well rounded.
\ep

\proof
We will assume as we may that all the orbits $Z_i$
have positive dimension.
When $I_\infty=\emptyset$, the well-roundedness of $B_T$
is a consequence of Example \ref{ex5}. Hence we will assume that
$I_\infty\neq\emptyset$.
When the height function $h_\infty$ is constant on $Z_\infty$,
the well-roundedness of $B_T$
is a consequence of Example \ref{ex1}.

When the height function $h_\infty$ is not constant on $Z_\infty$,
the well-roundedness of $B_T$ follows from Lemma \ref{lemex6} (4)
and of the linearity of the action of
each $G_i$ on $V_i$.
\hb\vs

Although we stated the above proposition
only for characteristic $0$ fields,
when all the $k_i$ have positive characteristic,
the height balls are also well-rounded by
Example \ref{ex5}.

\section{Applications}
\label{secapplication}

We will be applying the following theorem and corollary
to the above examples of well rounded sequences.
\bt
\label{thwellrounded}
Let $I$ be a finite set. For each $i\in I$, let $k_i$ be a local field
of characteristic not $2$,
$G_i$  the group
of $k_i$-points of a semisimple algebraic $k_i$-group,
$H_i\subset G_i$ the $k_i$-points of a symmetric $k_i$-subgroup.
Set $G_I:=\prod_{i\in I}G_i$, $H_I:=\prod_{i\in I}H_i$,
$Z_I:=H_I\ba G_I$ and $z_0=[H_I]$.
Let $\Ga$ be an irreducible lattice of $G_I$
such that $\Ga_H:=\Ga\cap H_I$ is a lattice in $H_I$. Then
for any  sequence $B_n$ of well-rounded subsets of $Z_I$
with volume tending to infinity, we have, as $n\to \infty$,
$$\#(z_0\Ga \cap B_n)\sim \tfrac{\mu_Y(Y)}{\mu_X(X)}\mu_{Z_I}(B_n),$$
where $X=\Ga \ba G_I$,  $Y= \G_H\ba H_I$
and the volumes are computed using invariant measures
as in (\ref{eqnin}).
\et

\proof
Use Corollary \ref{corequihorbit} and
Proposition \ref{prowellrounded} with $\ph_n:={\bf 1}_{B_n}$.
\hb\vs

In the product situation
of $Z_I=Z_{I_0}\times Z_{I_1}$, we will be taking a well-rounded
sequences of $Z_I$ which are products of a fixed compact subset
in one factor $Z_{I_1}$ and a well-rounded sequence of subsets in
 the other factor $Z_{I_0}$. This will give us equidistribution results
in the space $Z_{I_1}$ when $Z_{I_0}$ is non-compact.

\bc
\label{corwellrounded}
Keeping the same hypothesis as in Theorem \ref{thwellrounded},
let $I=I_0\sqcup I_1$ be a partition of $I$.
Letting $B_n$ be a well-rounded sequence of subsets of $Z_{I_1}$
with volume going to infinity, consider the following discrete
multisets $Z(n)$ of $Z_{I_0}$:
\begin{eqnarray*}
\label{eqnzpn}
Z(n):=\{ z\in Z_{I_0}\mid
(z,z')\in  z_0 \G \cap (Z_{I_0}\times B_n)
\quad \text{for some $z'\in Z_{I_1}$}\}.
\end{eqnarray*} Then, as $n\to \infty$,
the sets $Z(n)$
become equidistributed in $Z_{I_0}$ with respect to
a suitably normalized invariant measure.
In fact, for any $\ph\in C_c(Z_{I_0})$,
$$\lim_{n\ra\infty}
\frac{1}{\mu_{Z_{I_1}}(B_{n})}\sum_{z\in Z(n)}\ph(z)=
\tfrac{\mu_Y(Y)}{\mu_X(X)}
\int_{Z_{I_0}}\ph \; d\mu_{Z_{I_0}} .$$

In particular, $Z(n)$ is non-empty for all large $n$.
\ec

Multiset means that the points
of $Z(n)$ are counted with multiplicity
according to the cardinality of the fibers
of the projection  $z_0\Ga \cap (Z_{I_0}\times B_n)\ra Z(n)$.
Since $z_0\Ga $ is discrete and $B_n$ is relatively compact,
we note that these fibers are finite and that $Z(n)$ is discrete in
 $Z_{I_0}$.
\vs

\proof
It suffices to prove the claim for non-negative
functions $\ph\in C_c(Z_{I_0})$.
Define a sequence  of functions $\ph_n$ on $Z$
 by $$\ph_n(z,z'):=(\ph\otimes {\bf 1}_{B_n})(z, z')=
\ph (z){\bf 1}_{B_n}(z') \quad\text{for }
 (z, z')\in  Z_{I_0}\times  Z_{I_1}.$$
By Example \ref{ex1}, this sequence $\ph_n$ is well-rounded
and
$$\sum_{z\in Z(n)}\ph(z)
=\sum_{z\in z_0 \Ga }\ph_n(z) $$
since $Z(n)$ is a multiset.

By Corollary \ref{corequihorbit}, we can apply
Proposition \ref{prowellrounded}
to the sequence $\ph_n$ and obtain
\centerline{
\mbox{ }\hfill
$\D \lim_{n\ra\infty}
{\textstyle \frac{1}{\mu_{Z_{I_1}}(B_{n})}}\sum_{z\in Z(n)}\ph(z)
=
\lim_{n\ra\infty}
{\textstyle \frac{1}{\mu_{Z_{I_1}}(B_{n})}}
\sum_{z\in z_0\Gamma}\ph_n(z)
=
{\tfrac{\mu_Y(Y)}{\mu_X(X)}}
\int_{Z_{I_0}}\ph \; d\mu_{Z_{I_0}}\; .$
\hb}
\vs

\br
\label{remapplication}
\rm In Propositions \ref{prosumintegral} and  \ref{prowellrounded},
one can replace the hypothesis
``the $L^1$-norm of $\ph_n$ goes to infinity'' by the hypothesis that
``the support of $\ph_n$ is non-empty and goes to infinity''
i.e. for every compact $C$ of $Z$, $\ph_n|_C$ is null for all $n$ large.
The proof is exactly the same.
A similar remark applies to Theorem \ref{thintro},
\ref{thwellrounded}
and Corollary \ref{corwellrounded}.
This remark is useful for the non-empty spheres in Example \ref{ex4},
since it avoids to check that their volume goes to infinity
with the radius.
\er

\noi{\it Proof of Theorems \ref{thintro}, \ref{thoeq}, \ref{eqqq}
and Corollary \ref{corintro}.}
We are now ready to prove the non-effective statements
in the introduction.

Theorem \ref{thintro} is a consequence of
Theorem \ref{thwellrounded}
with $I=S$, $G_S=\bg_S$, $H_S={\bf H}_S$ and $\Ga=\Ga_S$.
The only thing we have to check is that $\Ga_S$ is an
irreducible lattice in ${\bf G}_S$.
This is the following classical lemma \ref{lemirr}.

Theorem \ref{thoeq} is a consequence of Corollary \ref{corwellrounded}
with $I=S$ and $I_0=S_0$.
Note that the projection
$z_0 \Ga \to Z_{I_0}$ is injective and hence
the multiset $Z(n)$ is a set.

Theorem \ref{eqqq} is an immediate consequence
of Corollary \ref{corequihorbit}.

Corollary \ref{corintro} is a consequence of Theorem \ref{thintro}
and Proposition \ref{proex6}.
\hb\vs

Let ${\bf G}$ be a connected semisimple group
defined over a global field $K$,
and let $S$ be a finite set of places of $K$ containing all
archimedean places $v$ such that ${\bf G}_v$ is non-compact.
Recall that these conditions assure that the
subgroup ${\bf G} _{\c O_S}$ is a lattice in $\bg_S:=\prod_{v\in S} \bg_v$.

\bl
\label{lemirr}
Let $\Ga_S$ be a subgroup of finite index in ${\bf G}_{\c O_S}$.
Suppose that $\bg$ is simply connected, almost $K$-simple
and that $\bg_S$ is non-compact.
Then $\G_S$ is an irreducible lattice in
${\bf G}_S$ (see Definition \ref{defirreducible}).
\el

\noi{\it Proof.}
Since $\bg$ is simply connected
and $\bg_S$ is non-compact, then
$\bg$ has the strong approximation property
with respect to $S$, that is,
the diagonal embedding of $\bg_K$ is dense in the
$S$-adeles $\bg_{\mathbb A_S}$, i.e., the adeles without $S$-component
 (see \cite[Th. 7.12]{PR} for characteristic
$0$ cases and \cite{P} for the positive characteristic case).
Since $\bg$ is $K$-simple, it follows  that
$\G_S$ is an irreducible lattice in
${\bf G}_S$ \cite[Cor. I.2.3.2 \& Th. II.6.8]{Ma2}.
\hb\vs

\bq
For the rest of this paper, we will transform the proofs explained
in the above chapters into  effective proofs.
For that we need to control precisely all the error terms appearing
in these proofs. There are mainly four  error terms to control.
The first three come from
the  mixing property,
the wave front property and
the approximation of $\mu_Y$ by a smooth function.
Their control will give
the equidistribution speed of the translates of $\mu_Y$.
The last error term comes from the
 well roundedness of the balls $B_T$.
We will dedicate one section to each of these terms.
\eq

\section{Effective mixing}
\label{seceffmix}
In this section,
we introduce notations which will be used
through the section \ref{seceffapp} and
we describe an effective version (Theorem \ref{progmo}) of the
mixing property based on the uniform decay of matrix coefficients.
\vs

We let $K$ be a number field, ${\bf G}$ a connected
simply connected almost $K$-simple group and
${\bf H}$ a $K$-subgroup of ${\bf G}$ with no non-trivial $K$-character.
Let
$S$ be a finite set of places of $K$ containing all
the infinite places $v$ such that
$\bg_v$ is non-compact.
 We write $S_\infty$ and
$S_f$ for the sets of infinite and finite places in $S$ respectively.
We assume that $\bg_S:=\prod_{v\in S}\bg_v$ is non-compact.
Let $\Gamma_S$ be a finite index
subgroup of $\bg (\c O_S)$. Note that $\bh_S\cap \Gamma_S$
is a lattice in $\bh_S$.

Set
$X_S:=\Gamma_S\ba {\bf G_S}$ and $Y_S=\G_S\cap \bh_S \ba \bh_S$.
 Let $\mu_{X_S}$ and
$\mu_{Y_S}$ denote the invariant probability measures on $X_S$ and $Y_S$
respectively.
Set $Z_S:=\bh_S\ba {\bf G_S}$. For each $v\in S$, choose
 an invariant measure $\mu_{Z_v}$ on ${\bf H}_{v}\ba{\bf G}_v$
so that the invariant measure
 $\mu_{Z_S}:=\prod_{v\in S}\mu_{Z_v}$
 on $Z_S$ satisfies
$\mu_{X_S}=\mu_{Y_S}\,\mu_{Z_S}$ locally. For $S_0\subset S$,
we set $\mu_{Z_{S_0}}:=\prod_{v\in S_0}\mu_{Z_v}$.

By a smooth function on $X_S$
we mean a function which is smooth on each ${\bf G}_\infty$-orbit
and which is invariant under a compact open subgroup of ${\bf G}_f$.
The notation $C^\infty_c(\Gamma_S\ba {\bf G}_S)$
denotes the set of smooth functions with compact support on ${\bf G}_S$.

For each $v\in S$,
recall the ``Cartan'' decomposition due to Bruhat and Tits in
\cite{BT1} and \cite{BT2}:
one has  ${\bf G}_v=M_v
\Omega_vB_v^+M_v$ where $M_v$ is a good maximal
compact subgroup, $B_v^+$ a positive Weyl chamber
of a maximal $K_v$-split torus and $\Omega_v$ is a finite subset
in the centralizer of $B_v$.

For simplicity, we set ${\bf G}_{\infty}={\bf G}_{S_\infty}$
and ${\bf G}_f={\bf G}_{S_f}$.
We also set $M_\infty:=\prod_{v\in S_\infty} M_v$
and $M_f:=\prod_{v\in S_f} M_v$.

Let $X_1, \cdots, X_d$ be an orthonormal basis of the Lie algebra
of $M_\infty$
 with respect to an $\op{Ad}$-invariant
scalar product. We denote by $\c D$
the elliptic operator $\c D:=1-\sum_{i=1}^d X_i^2$.

Fix any closed embedding of ${\bf Z}={\bf H}\ba {\bf G}$ into
a finite dimensional vector space ${\bf V}$ defined over $K$;
such an embedding always exists by the well known theorem
of Chevalley.
To measure how far an element $z\in {\bf Z}_S$ is from
the base point $z_0=[\bh_S]$ in $Z_S:=\bh_S\ba \bg_S$,
we may use a height function
  \begin{equation}\label{hedd}\H_S(z):=
  \prod_{v\in S}\|z_v  \|_v\end{equation}
where
$\|\cdot \|_v$ is a norm
on ${\bf V}_{v}$.
This norm is assumed to be euclidean when $v$ is an infinite place
and a max norm when $v$ is a finite place.
Note that the height function $\H_S:Z_S \to \mathbb R^+$ is
 a proper function.

\bt
\label{progmo} There exists $\kappa>0$ and $m\in \m N$ such that
for
any  open compact subgroup $U_f$ of $\bg_f$,
there exists $C_{U_f}>0$ satisfying that
 for any $\psi_1, \psi_2\in C_c^\infty
(X_S)^{U_f}$ and any $g\in \bg_S$,
\begin{equation*}
\left|\langle g \psi_1, \psi_2\rangle -\int_{X_S} \psi_1 d\mu_{X_S} \;
\int_{X_S} \psi_2d\mu_{X_S}\right| \le
C_{U_f}  \H_S(z_0 g )^{-\kappa}
 \| \mathcal D^m(\psi_1)\|_{L^2}
\|\mathcal D^m(\psi_2)\|_{L^2} .\end{equation*}
\et

\proof
The above claim is a straightforward consequence
of Theorem 2.20 of \cite{GMO} based on the results
of \cite{Clo} and \cite{Oh1}.
\cite[Theorem 2.20]{GMO} relies on the following hypothesis:
''the only character appearing in $L^2(\G_S\ba \bg_S)$
is the trivial one''. This hypothesis is satisfied here since
the non compactness of $\bg_S$ and the
simply-connectedness of ${\bf G}$ imply
the irreducibility of $\Ga_S$ by Lemma \ref{lemirr}.

The conclusion of  \cite[Theorem 2.20]{GMO} is the above claim where
$ \H_S^{-\kappa}$ is replaced by a function
$\wt\xi_{\bf G}$
which is a product over $v\in S$
of  bi-$M_v$-invariant functions $\xi_v'$
satisfying
$$\xi'_v(a)\le \prod_{\alpha\in Q_v}\alpha(a)^{-1/2+\e}\quad
\text{for all $a\in B^+_v$} $$
where $Q_v$ is a maximal strongly orthogonal system of the root
system of $({\bf G}_v, B_v)$

We only have to check that this function
$\wt\xi_{\bf G}$ is bounded by a multiple of $ \H_S^{-\kappa}$.
For that, denote by $\rho$ the representation
of ${\bf G}$ into ${\bf GL}({\bf V})$
such that the stabilizer of $z_0\in {\bf V}_K$ is ${\bf H}$ and
choose a weight $\lambda$ larger on $B_v^+$ than any weight of $\rho$.
Then there exists a positive integer k such that, for all $a\in B_v^+$,
 $$\|z_0\rho( a) \|_v\le \|z_0\|_v \, \|\rho(a)\|_v
\le \|z_0\|_v \, |\lambda (a)|_v \le \|z_0\|_v \,
\prod_{\alpha\in Q_v}\alpha(a)^k
.$$
Since $M_v$ and $\Omega_v$ are compact subsets,
by the continuity, this implies that there exists $\kappa >0$ and $c>0$
 such that
$$\xi_v'(g)\le  c\, \|z_0 \rho(g)\|_v^{-\kappa}$$
for all $g\in {\bf G}_v$. This implies our claim.
\hb

\section{Injective radius and the approximation
by smooth functions}
\label{seceffnei}

The aim of this section is to get an effective upper bound
on the volume of the set of points in $Y_S$ with small injectivity radius
in $X_S$ and approximate the characteristic function
\vs

Fix a closed embedding $\bg \hookrightarrow {\bf GL}_N$.
We
may consider each element $g$ of ${\bf G}_S$ as an $|S|$-tuples  of
$N\times N$ matrices $g_v$. We also fix a norm $\|.\|_v$
on each of these $K_v$-vector spaces of matrices.

For $x\in X_S$,
consider the projection map $p_x:\bg_S \to X_S$ given by $g\mapsto xg$.
The injectivity radius $r_x$
is defined to be
$$r_x= \sup\{r>0: p_x|_{B_r \times M_f}
\text{ is injective}\}$$
where $\D B_r=\{g\in \bg_\infty : \max_{v\in S_\infty}\|g_v-e\|_v \leq r\; \}$.

Of course, this definition makes sense only when $S_\infty$ is
non-empty. This does not matter since, when $S_\infty$ is empty,
$X_S$ is compact.

\begin{Lem}
\label{lemrq} Suppose $S_\infty\neq\emptyset$.
For any  $x\in X_S$, one has
 $r_x>0$.
\end{Lem}

\proof
Since $\Ga_S$ does not meet ${\bf G}_f$
and ${\bf G}_f$ is normal in ${\bf G}_S$,
the group ${\bf G}_f$ acts freely on $X_S$.
Hence $p_z|_{\{ e\}\times M_f}$ is injective.
Since $M_f$ is compact  and $p_x$ is locally injective,
$p_x|_{B_r \times M_f}$ is still injective for some small $r>0$.
\hb\vs

Moreover we have a quantitative version of the above lemma.

\begin{Lem}
\label{inject}
Suppose $S_\infty\neq\emptyset$.
There exist
  $c_1>0$, $p_1>0$ such that for all
sufficiently small $\e>0$,
$\mu_{Y_S}(\{ y\in Y_S\mid r_y<\eps\})\le c_1 \e^{p_1} .$
\end{Lem}

\proof
We use the reduction theory for $\bh_S$ (cf. \cite{PR}). We first recall
what a Siegel set is. Let ${\bf A}$ be a maximal $K$-split torus of
$\bh$ and $\bf P$ a minimal parabolic subgroup containing $\bf A$.
Then ${\bf P}={\bf N}{\bf R}{\bf A}$ where ${\bf R}$ is a
$\q$-anisotropic reductive subgroup and ${\bf N}$ the unipotent
radical of ${\bf P}$. Set ${\bf A}_\infty:=\prod_{v\in S_\infty}
{\bf A}(K_v)$ and similarly ${\bf N}_\infty$ and ${\bf R}_\infty$.
Denoting by $\Delta$ the system of simple roots of $\bh_\infty$
determined by the choice of ${\bf P}$, we set for $t>0$,
$$A_t=\{a\in {\bf A}_\infty: \alpha(a)\ge t
\quad\text{for all $\alpha\in \Delta$}\}.$$ Then for a compact
subset $\omega\subset {\bf N}_\infty {\bf R}_\infty$
and a
maximal compact subgroup $K_0$ of $\bh_S$,
the set
$\Sigma_t:=\omega A_t K_0$ is called a Siegel set.
Now the reduction theory says
that there exist $h_1, \cdots, h_r\in \bh_S$, and a Siegel set
$\Sigma_{t_0}= \omega A_{t_0} K_0$ such that
$$\bh_S=\cup_{i=1}^{r} (\bh_S\cap \G_S) h_i
 \Sigma_{t_0}. $$

Let $M'_f:=\cup_{i=1}^r h_i M_f h_i^{-1}$.
As in Lemma \ref{lemrq},
there exists $\e_0>0$ such that
\begin{equation}\label{of}
\Gamma_S\cap B_{\e_0}M'_f=\{e\}. \end{equation}

Set
$$C_\e:=\cup_{i=1}^r h_i \{ wak\in \Sigma_{0}
: t_0\le \alpha(a) \le
\e^{-r_0} \text{ for each $\alpha\in \Delta$ }\} ,$$
where $r_0>0$ is chosen independent of $\eps$,
so that, for all $g$ in $C_\eps$ and $v\in S_\infty$, one has
$$
\| g_v\|_v\leq \eps^{-1/4}
\;\;{\rm and}\;\;
\| g^{-1}_v\|_v\leq \eps^{-1/4}\; .
$$
Let $Y'_\e $ denote the image of $C_{\e}$
in $Y_S$ under the projection $\bh_S\to Y_S$.
The integration formula \cite[p. 213]{PR}
shows that for some constant $c_1>0$ and $p_1>0$,
$$\mu_{Y_S}(Y_S\! -\! Y'_\e)\le c_1 \e ^{p_1} .$$
Hence it is enough to show that for all $z\in Y'_\e$,
one has $r_z\geq\eps$.
Suppose $p_{z}(x)=p_{z}(y)$ with $x, y\in B_\e M_f$ and write
$z=\G_S g$ for some $g\in C_{\e}$.
We want to prove that $x=y$.

The element  $\gamma:= g x y^{-1}g^{-1}$ belongs to $ \G_S$.
Moreover, for some fixed constant $c>1$, one has, for all $v\in S_\infty$,
\begin{align*}
\|\gamma-e\|_v &=
\|g_v (x_v -y_v)y_v^{-1} g_v^{-1}\|_v
\\
&\le  c \,\e ^{-1/2}\, \|x_v-y_v\|_v\; \|y^{-1}_v\|_v
 \\ &\le  c^{2}  \e^{1/2}.
\end{align*}
But  the finite component of $\ga$ is in $M'_f$,
hence, $\ga$ is in $B_{c^2\eps^{1/2}}\times M_f$ and
one gets from (\ref{of})
that, for  $\e <c^{-4}\e_0^2$, one has
$\gamma=e .$
Therefore $x=y$ as well.
\hb\vs

For all $v\in S$, we choose
a small neighborhood $\g s_v$ of $0$ in
a supplementary subspace of the Lie algebra $\g h_v$  in  $\g g_v$ and
set $\g s:=\prod_{v\in S}\g s_v$.
The set  $W:=exp( \g s)$ is then a transversal to $\bh_S$
in $\bg_S$.
We set $\mu_W$  the measure on $W$ such that $d\mu_{X_S}=d\mu_{Y_S} d\mu_W$
locally.

Recall that $B_\e$ denotes the ball of center $e$
and radius $\e$ in ${\bf G}_\infty$
and let $U_\e$ be the ball of center $e$ and radius $\e$ in ${\bf G}_f$:
\begin{eqnarray}\label{bed}
B_\e
&=&
\{ g\in \bg_\infty  \mid \max _{s\in S_{\infty}} \| g-e\|_v \leq \eps \},\\
U_\e
&=&
\{ g\in \bg_f  \mid \max _{s\in S_f} \| g-e\|_v \leq \eps \}.
\end{eqnarray}

We fix  $\e_0>0$ small. For $\e$ small we let
$$H_\e:=\bh_S\cap B_\e U_{\e_0}\quad\text{and}\quad
W_\e:=W\cap B_\e U_{\e_0}$$
so that the multiplication $H_\e\times W_\e\ra H_\e W_\e$ is an homeomorphism
onto a neighborhood of $e$.
Fix $m > \dim \bg_\infty$ and $\kappa>0$ satisfying Theorem \ref{progmo} and
fix $l\in \m N$ as in Lemma \ref{wave}.
We can assume that $U_{l\e_0} \subset U_f\cap M_f$
and let $U'_f$ be an open subgroup of $\bg_f$
such that $U_{\e_0}U'_f=U_{\e_0}$.

\bl
\label{lemtra}
Let  $Y_\eps:=\{y\in Y_S \mid \mbox{\rm the map}\; g\mapsto yg \;
\mbox{\rm is injective on}\; H_\e W_\e\}$.
There exist $c_1>0$ and $ p_1>0$ such that for all small $\e>0$,
$$\mu_{Y_S}(Y_S\!-\! Y_\eps)\le c_1 \e^{p_1} .$$
\el

\proof
When $S_\infty\ne \emptyset$,
since
$B_{\eps}B_{\eps}\subset B_{\e^{1/2}}$ for $\e$ small,
the set $Y_\eps$ contains  the set of points $y$
such that $r_y\geq \eps^{1/2}$.
Just apply then Lemma \ref{inject}.

When $S_\infty=\emptyset$, $X_S$ is compact,
hence $Y_\eps$ is equal to $Y_S$ for $\e_0$ and $\e$ small.
\hb\vs

The following proposition provides the approximation of
the characteristic function ${\bf 1}_{Y_\e W}$
by a smooth function  $\ph_\e$ with the controlled Sobolev norm.

We first recall
the Sobolev norm ${\c S}_m(\psi)$ of a function $\psi\in C_c^\infty(X_S)$.
Choose a basis $X_1,...,X_n$
of the Lie algebra of ${\bf G}_\infty$.
For each $k$-tuple of integers $a:=(a_1,...,a_k)$  with $1\leq a_i\leq n$,
the product $X_a:=X_{a_1}\ldots X_{a_k}$
defines a left-invariant differential operator
on ${\bf G}_S$, hence a differential operator on $X_S$.
By definition ${\c S_m}(\psi)^2=\sum_a\|X_a\psi\|_{L^2}^2$
where the sum is over all the $k$-tuples $a$ with $0\leq k\leq m$
and where $X_{_\emptyset}\psi$ stands for $\psi$.
\bp
\label{lemphe}
There exist $p_2>0$
such that, for all sufficiently small $\e$, one can choose
\\
- a non-negative smooth function $\rh_\e$ on $W$
with support in $W_\e$ such that $\int_W\rho_\e=1$, \\
- a non-negative smooth function $\tau_\e$ on $Y_S$
with support in $Y_\e$ such that $\tau_\e \le 1$ on $Y_S$
and $\tau_\e|_{Y_{4\e}} =1$.
\\
- Moreover, let $\ph_\e$ be the function on $X_S$ defined by
$$
\ph_\e(x)=\sum_{\{(y,w)\in Y_\e\times W_\e \mid yw=x\}}
\tau_\e (y)\rho_\e(w).
$$
The choices can be made so that
$\ph_\e$ is $U'_f$-invariant and
$\c S_m(\ph_\e)\leq \e^{-p_2}$.
\ep

We remark that the sum defining $\ph_\e$ is a finite sum and hence
$\ph_\e$ is well defined.

To prove Proposition \ref{lemphe},
we first need a lemma which constructs some test functions
$\al_\e$ near $e$.

\begin{Lem}
\label{lemsobolev}
For a given $m\geq 0$, there exists $p\in \m N$, such that,
for all sufficiently small $\e>0$, one can choose
smooth non-negative functions $\be_\e$ on $H_\e$
and smooth non-negative functions $\rho_\e$ on $W_\e$
satisfying the following: \\
- one has $\be_\e\geq 1$ on $H_{\e^2}$.\\
- one has $\int_W \rho_\e d\mu_W =1$\\
- if
$\al_\e$ denotes the smooth function on $H_\e W_\e$ given by
$\al_\e (hw)=\be_\e(h)\rho_\e(w)$, then
  $\al_\e$ is $U'_f$-invariant
and $ \c S_m(\al_\e)\le \e^{-p}$.
\end{Lem}

\proof
The general case reduces to the case
of $S=S_\infty$, by considering tensor products
with  characteristic functions
of $U_{\e_0}\cap H$ and of $U_{\e_0}\cap W$.
Hence, we can assume that $S=S_\infty$ so that
$\bg_S$ is a real Lie group. Set $d= \dim W$.
Fix some smooth non-negative functions $\be$ on
$\g h:=\oplus_{v\in S}\g h_v$ and $\rho$ on $\g s$
with support in a
sufficiently small neighborhood of $0$
such that $\be (0) >1$ and $\int_{\g s}\rho =1$.
Then, for suitable constants $c_\e>0$ converging to $1$,
the functions given by
$$
\be_\e(exp( X))=\be (X/\e)
\;\; \text{ and } \;\;
\rho_\e( exp( Y))=c_\e \eps^{-d}\rho (Y/\e),
$$
for $X$ (resp. $Y$) in
a fixed compact neighborhood of $0$ in $\g h$ (resp. $\g s$),
satisfy the properties listed above.
\hb\vs

\noi {\it Proof of Proposition \ref{lemphe}.}
We choose the function $\rho_\e$ from Lemma \ref{lemsobolev}.
To construct the function $\tau_\e$,
consider
a maximal family $\c G_\e$ of points
$y\in Y_\e$ such that the subsets
$y H_{\e^3}$ of $Y_S$ are disjoint and meet $Y_{2\e}$
and let $\c F_\e\subset \c G_\e$ the subfamily for which
$yH_{\e^3}$ meets $Y_{4\e}$.
For all $y\in \c G_\e$  the volumes
$\mu_{Y_S} (y H_{\e^3})$ are equal and
of order $\e^{3d}$ with $d=\text{dim}(\bh_\infty)$.
Since $\mu_{Y_S}(Y_S)=1$,
the cardinality of $\c G_\e$ is at most $O(\e^{-3d})$.

For $y\in \c G_{\e}$ we define a test function
$\be_{y,\e}$ on $Y_S$ with support on $yH_\e$ by
$\be_{y,\e}(yh)=\be_\e(h)$ and
let
$\be_{\c G,\e}:=\sum_{y\in \c G_\e}\be_{y,\e}$.
Since $B_{\e^3}B_{\e^3}\subset B_{\e^2}$,
the sets
$y H_{\e^2}$, $y\in \c G_\e$,  cover  $Y_{2\e}$. Hence
$\be_{\c G,\e}\ge 1$ on $Y_{2\eps}$.

For each $y\in \c F_{\e}$, consider the
test function $\tau_{y,\e}$ on $Y_S$ with support in $Y_{2\eps}$
given by
$\tau_{y,\e}:=\be_{y,\e}/\be_{\c G,\e}$ on $Y_{2\eps}$
and set
$$\textstyle \tau_\e :=\sum_{y\in \c F_\e}\tau_{y,\e}\, .$$
Note that $0\le \tau_\e \le 1$
 on $Y_S$, $\tau_\e|_{Y_4\e}=1$ and $\tau_\e|_{Y_S-Y_\e}= 0$.
 For $y\in \c F_{\e}$,
we also define
the test function $\ph_{y,\e}$  on $X_S$
with support on $yH_\e W_\e$ given by
$$\ph_{y,\e}(yhw):=\tau_{y,\e}(yh)\rho_\e(w)=\al_\e(hw)/\be_{\c G,\e}(yh).$$
These functions $\ph_{y,\e}$ are
well-defined since $y$ belongs to the set $Y_\e$ given
by Lemma \ref{lemtra}.
By construction, we have
$$\textstyle \ph_\e =\sum_{y\in \c F_\e}\ph_{y,\e}\, .$$
It follows from $\c S_m(\al_\e)\leq \e^{-p}$ that there exists $p_0>0$
such that
$$\D\max_{y\in \c F_\e}\c S_m(\ph_{y,\e}) =O(\e^{-p_0}) $$
and hence $\c S_m(\ph_{\e}) =O(\e^{-3d-p_0}) $.
\hb

\section{Effective equidistribution of translates of $\bh_S$-orbits}
\label{eff}

The goal of this section is to prove Theorem \ref{transff},
or its stronger version Theorem \ref{thmeet} below.
 This is an
effective version of Proposition \ref{proequihorbit}
on the equidistribution of translates of $\bh_S$-orbits in $X_S$.
\vs

\begin{Def}
\label{dthmeet} We say that the translates $Y_Sg$
are effectively equidistributed in $X_S$ as $g\to \infty$ in $Z_S$ if
there exists $m\in \m N$ and $r>0$ such that, for any
compact open subgroup $U_f$ of $\bg_f$ and
any compact subset $C$ of $X_S$,
there exists  $c=c(U_f,C)>0$ satisfying that
for any smooth function $\psi\in C^\infty_c(X_S)^{U_f}$
with support in $C$, one has for
all $g\in \bg_S$
\begin{eqnarray}
\label{eqneet}
|\int_{ Y_S} \psi (yg)\,
d\mu_{Y_S}(y)-\int_{X_S}\psi\,d\mu_{X_S}|
\le c \,S_m(\psi)\, \H_S( z_0 g)^{-r} .
\end{eqnarray}
\end{Def}

Assume further that
$\bh$ is a symmetric $K$-subgroup of $\bg$.
Taking the product of  the polar
decompositions  $\bg_v=\bh_v A_v K_v$ given in
Lemma \ref{lemkah} over $v\in S$,
we obtain a polar decomposition of the shape
$$\bg_S=\bh_S A_S K_S \;\;\text{ and  we set }\;\; F_S:=A_SK_S.$$

The following effective version of the wavefront property
\ref{defwavefront} is a main technical reason why our proof
of Theorem \ref{thmeet} works for $\bh$ a symmetric subgroup.
\begin{Lem}
\label{wave}
There exists $l\in \m N$ such that
 for all small $\e,\e'>0$ and
 all $g\in F_S$,
\begin{equation}
\bh_S  \;B_{\e/l} U_{\e'/l}\; g\subset \bh_S \, g \,B_\e  U_{\e'}\;.
\end{equation}
\end{Lem}
\proof  We only have to check this separately at each place $v$.
This statement is then a strengthening of
Proposition \ref{prowavefront}  on
the wavefront property and is
an output of the proof of this Proposition.
\hb

\begin{Thm}\label{thmeet} If ${\bh}$ is a symmetric $K$-subgroup
of $\bg$,
then the translates $Y_Sg$
are effectively equidistributed in $X_S$ as $g\to \infty$ in $Z_S$.
\end{Thm}
\proof
Since $\bg_S=\bh_S F_S$,
 it suffices to prove the above claim for $g\in F_S$.
We may also assume that $\int_{X_S} \psi\, d\mu_{X_S} =1$.
We want to bound $|I_g-1|$ where
$$I_g:=\int_{Y_S}\psi (yg)d\mu_{Y_S}(y) .$$
We  follow the proof of Proposition \ref{proequihorbit}.
The main modification will be to replace the characteristic function
${\bf 1}_{Y_\e W}$ by the test function $\ph_\e$
constructed in  Proposition \ref{lemphe}.
By the same argument as in section \ref{sechorbits},
but using  the stronger version \ref{wave} of the wavefront lemma,
we have that for all small $\e>0$
and for any $w\in W_\e$
\begin{equation}
\label{eqnig3}
|I_g-\int_{Y_S}\psi (ywg)d\mu_{Y_S}(y) | \le
 l \e \, C_\psi .
\end{equation}
Here $C_\psi$ is the Lipschitz constant at $\infty$, i.e. the smallest
constant such that
for all $\e >0$,
$|\psi(xu)-\psi(x)|\le C_\psi \,\e\quad\text{for all $x\in X_S$ and
$u\in B_\e$.}$

Set $\tau_\e$, $\rho_\e$ and $\ph_\e$ the functions
constructed in  Proposition \ref{lemphe} and$$
J_{g,\e}:=\int_{W_{\e}} \int_{Y_S} \psi(ywg) \rho_\e(w)
d\mu_{Y_S}(y) d\mu_W(w).$$

By integrating \eqref{eqnig3} against $\rho_\e$, we obtain
$$|I_g - J_{g,\e}|\le  l \e \, C_\psi .$$

Set also
\begin{eqnarray*}
K_{g,\e}&:=&
\int_{X_S} \psi (xg)
\ph_\e(x)\,
 d\mu_{X_S}(x)\\
&=&
\int_{W_{\e}} \int_{Y_S} \psi(ywg) \tau_\e(y)\rho_\e(w)
d\mu_{Y_S}(y)\, d\mu_W(w).
\end{eqnarray*}

Noting that
$ \tau_\e(y)=1$ for $y\in Y_{4\e}$,
we have for some $c_1,p_1>0$
\begin{eqnarray*} |J_{g,\e}-K_{g,\e}|&=&
\left| \int_{W_\e } \int_{Y_S}\psi(ywg)
(1-\tau_{\e}(y)) \rho_\e(w) d\mu_{Y_S}(y)\,d\mu_W(w) \right| \\
 &\leq &
 2 \mu_{Y_S}(Y_S\! -\! Y_{4\e}) (\int_{W_\e}\rho_\e )
\|\psi\|_\infty\,
  \\
& \le& c_1\, \e^{p_1} \|\psi\|_\infty\; .
\end{eqnarray*}
Note that
$K_{g,\e}=\langle g.\psi, \ph_\e \rangle .$
Since $\ph_\e$ and $\psi$ are $U'_f$-invariant,
by Theorem \ref{progmo} and
Proposition \ref{lemphe} ,
we deduce for some $c',c_2,p_2>0$
\begin{eqnarray*}
|K_{g,\e} - \int_{X_S}
\ph_\e d\mu_{X_S}|
 &\le&
 c'\,\| {\c D}^m(\psi)\|_{L^2} \|{\c D}^m(\ph_\e)\|_{L^2}
\H_S(z_0g)^{-\kappa}\\
&\leq &
  c_2\e^{-p_2} \|{\c D}^m(\psi)\|_{L^2}
\H_S(z_0g)^{-\kappa} .
\end{eqnarray*}

Moreover, one has
$$\left|\int_{X_S}
\ph_\e d\mu_{X_S} -1\right|=\left|
\int_{Y_S} \tau_{\e}(y)\, d\mu_{Y_S}(y)-1
\right|
\leq \mu_{Y_S}(Y\! -\! Y_{4\eps})
\leq c_1\, \eps^{p_1}\, .$$

Since  $C$ is compact,
the $C^1$-norm of a $U'_f$-invariant function $\psi$ supported on $C$
is bounded  above by a uniform multiple of a suitable Sobolev norm
as in \cite[Theorem 2.20]{Au} i.e., one has
an inequality
$$
\max(\|\psi\|_\infty,C_\psi)\le c''\, \c S_{2m}(\psi)
$$
with $c''=c''(U_f,C)>0$ independent of $\psi$. Hence,
putting all these upper bounds together and using the inequality
$1\leq \|\psi\|_\infty$,
we get
\begin{eqnarray*}|I_g -1|
&\le & \textstyle
|I_g-J_{g,\e}|+|J_{g,\e}-K_{g,\e}|+|K_{g,\e}-\int\ph_\e|+
|\int\ph_\e -1|\\
&\le &
l \e \, C_\psi + c_1 \e^{p_1} \|\psi\|_\infty+
c_2  \e^{-p_2}  \|\c D^m(\psi)\|_{L^2} \,
\H_S(z_0g)^{-\kappa} +c_1\e^{p_1} \\
&\le &
(c_1'\e^{p_1}   +
c_2'\e^{-p_2}  \H_S(z_0g)^{-\kappa})\mathcal S_{2m}(\psi).
\end{eqnarray*}
Note in the above that the positive constants $ c_i',p_i$, $i=1,2$,
 are independent of $\psi$.

Now by taking $\e=\H_S(z_0g)^{-r/p_1}$ with
 $r=\frac{\kappa \, p_1}{p_1+p_2}$, we obtain as required
$$|I_g-1|\le c \, \mathcal S_{2m}(\psi)\;\H_S(z_0g)^{-r}\, .$$
This concludes the proof.
\hb
\vs

\noi {\bf Remarks}
\begin{enumerate}\item One could also,
as an output of our proof, compute
explicitly $m$ and $r$ and describe  how
the constant $c$ depends on the compact sets $U_f$ and $C$.
\item
Note that the above theorem \ref{thmeet}
is precisely the effective version
of Proposition \ref{proequihorbit},
since we have shown that the effective mixing theorem
\ref{progmo}
together with the effective wave front
lemma \ref{wave} imply the effective equidistribution
of $Y_S g$.
\end{enumerate}

\section{Effective counting and equidistribution}
\label{sf}

The following definition is an effective version of Definition
\ref{defwellrounded}. Recall that
$B_\e=B(e,\e)$ is the ball of center $e$ and
radius $\e$ in ${\bf G}_{\infty}$ \eqref{bed}
and that $\H_S$ is a height function on $Z_S$ as defined in \eqref{hedd}.
\bd
\label{defeffwel}
A sequence of subsets $B_n$ in $Z_{S}$
is said to be effectively well-rounded if
\begin{enumerate}
\item it is invariant under a compact
open subgroup of $\bg_{ f}$,
\item there exists $\kappa>0$
such that, uniformly for all $n\geq 1$ and all $\e\in ]0,1[$,
\begin{equation*}
\label{bn}
\mu_{Z_S}(B_{n,\e}^+- B_{n, \e}^-)=
O(\e^\kappa \mu_{Z_S}(B_n))\end{equation*}
where $B_{n, \e}^+= B_nB_\e$ and $B_{n, \e}^-=\cap_{u\in B_\e} B_n u$,
\item
 for any $k>0$, there exists $\delta>0$
such that, uniformly for all $n\gg 1$
 and all $\e\in ]0,1[$,
 one has
\begin{equation*}
\label{bnhhh}
\int_{B_{n,\eps}^+} \H_{S}^{-k}(z)
\; d\mu_{Z_{S}}(z) =O(\mu_{Z_S}(B_{n})^{1-\delta}) .\end{equation*}
 \end{enumerate}
\ed
If $S_\infty$ is empty, then the assumption (2) is void.
\vs

A subset $\Om$ of $Z_S$ is said to be {\it effectively well-rounded}
if the constant sequence $B_n=\Om$ is effectively well-rounded.
This means that $\Om$ is of non-empty interior and that
the volume $\mu_{Z_S}(\partial_\eps\Om)$
of the $\eps$-neighborhood of the boundary of $\Om$
is a $O(\eps ^\kappa)$ for $\eps$ small.
For instance, a compact subset of $Z_{S_\infty}$
with piecewise smooth (or even piecewise ${\rm C}^1$) boundary is
effectively well-rounded in $Z_{S_\infty}$.

\begin{Thm}
\label{thspeed}
Suppose that the translates $Y_Sg$ become effectively
equidistributed in $X_S$ as $g\to\infty$ in $Z_S$.
Then for any effectively well-rounded sequence of subsets $B_n$ in
$Z_{S}$ such that $\vol(B_n)\ra\infty$ there exists a constant
$\delta_0>0$ such that
\begin{equation*}\label{omega}\# z_0\Gamma_S\cap  B_n =
\Vol(B_n)(1+O(\Vol(B_n)^{-\delta_0})) .\end{equation*}
\end{Thm}

\proof
Set $\Ga_{\! S,H}:= \G_S\cap \bh_S$. As in sections
\ref{secsumintegral} and \ref{seccounting}, we define
a function $F_n$ on $X_S=\G_S\ba \bg_S$ by
 $$F_n(x_0g)=\sum_{\gamma\in \Ga_{\! S,H}\ba \G_S}
{\bf 1}_{B_n}(z_0\ga g)
\;\; ,\; {\rm for} \; g\in \bg_{S}
.$$
For instance, one has
$$F_n(x_0)=\# z_0\Gamma_S\cap B_n .$$
Let $m$ and $r$ be the integers given by Theorem \ref{thmeet}
and $U_f$ a compact open subgroup of $\bg_{f}$.
By Lemma  \ref{lemsobolev},
there exists $p>0$, a smooth $U_f$-invariant
function $\al_\e$  on $\bg_S$,
supported on
$B_\e U_f$
such that $\int_{\bg_S} \al_\e=1$ and
$\mathcal S_m (\al_\e)\le \e ^{-p}$.
Here we take $\e$ and
$U_f$ small enough so that $B_\e  U_f$ injects to $X_S$, and hence
we may consider $\al_\e$ as a function on $X_S$.

We also introduce the functions $F_n^{\pm}$ on $X_S$:
\begin{eqnarray*}
\label{eqnfnpmeff}
F_{n,\eps}^\pm(x_0g)=\sum_{\ga\in \Ga_{S,H}\ba \Ga_S}
{\bf 1}_{B_{n,\eps}^\pm}( z_0 \ga g)\;
\;\; ,\; {\rm for} \; g\in \bg_{S} .
\end{eqnarray*}
Then
$$ F_{n, \e}^- (x_0g) \le F_n(x_0)\le
F_{n, \e}^+(x_0g)\quad\text{
for all $g\in B_\e\times U_f$}$$ and hence
$$\langle F_{n, \e}^-, \al_\e \rangle\le F_n(x_0)\le \langle
F_{n, \e}^+, \al_\e \rangle .$$
Note that $$\langle F_{n, \e}^\pm, \al_\e\rangle =
\int_{B_{n,\eps}^\pm}\left(
\int_{Y_S} \al_\e (yg) d\mu_{Y_S}(y)\right)d\mu_{Z_S} (z_0g) .$$

Set $v_n:=\vol(B_n)$ and $v_{n,\e}^\pm:=\vol(B_{n,\e}^\pm )$.
Then by Theorem \ref{thmeet} and the assumptions
(2) and (3) of the definition
 \ref{defeffwel},
there exist positive constants $\kappa$, $\delta$ and  $c_i$
such that for all $n\gg 1$ and small $\eps>0$,
\begin{align*}
|\langle F_{n, \e}^\pm, \al_\e\rangle  -v_n|\le &
\int_{B_{n,\eps}^+ }\left|\int_{Y_S} (\al_\e  (yg) -1)
 d\mu_{Y_S}(y)\right| d\mu_{Z_S} (z_0g) +
(v_{n,\e}^+-v_{n,\e}^-)\\
\le &  c\,\c S_m(\al_\e)
\int_{B_{n,\e}^+ }\H_S (z)^{-k}d\mu_{Z_S} (z) +
c_2 \e^\kappa v_n \\
\le &  c_1 \e^{-p} v_n^{1-\delta}
+c_2 \e^\kappa v_n.
\end{align*}
Setting $\de_0:=\frac{\de}{1+p/\kappa}$ and
choosing $\e=v_n^{-\de_0/\kappa}$
we get $F_n(x_0)=v_n(1+O(v_n^{-\delta_0}))$.
\hb

\begin{Cor}
\label{corspeed} Let $S=S_0\sqcup S_1$ be a partition of
$S$. There exist $\de_0,c >0$ such that
for any effectively well-rounded sequence of subsets $B_n$ in $Z_{S_1}$
whose volumes $v_n:=\mu_{Z_{S_1}}(B_n)$ tend to $\infty$ and
for any compact effectively well-rounded subset $\Om$ of $Z_{S_0}$,
we have
\begin{equation*}
\# z_0\Gamma_S\cap (\Omega\times B_n) =
 \tfrac{\mu_{Y_S}(Y_S)}{\mu_{X_S}(X_S)}\;
v_n\,\mu_{Z_{S_0}}(\Omega)\,(1+O(v_n^{-\delta_0})) \, . \\
\end{equation*}
\end{Cor}

\proof To apply Theorem \ref{thspeed},
we only have to check that the sequence $A_n:=\Om\times B_n$
of subsets of $Z_S$ is effectively well-rounded, which is straightforward.
\hb

\section{Effective well-roundedness}
\label{seceffwel}

In this section, we  give explicit examples of
effectively well-rounded families  (see Definition
\ref{defeffwel}).

We keep the notations for $K$, $S$, $\bg$, $\bh$, $Z$,
$\H_S(z)=\prod_{v\in S} \|z \|_v$ etc.,
 from the beginning of section \ref{seceffmix}.

We also set for $T >0$,
\begin{equation}\label{sphered}
B_v(T):=\{ z\in Z_v\mid \| z\|_v\leq T\}
\quad\text{and}\quad
S_v(T):=\{ z\in Z_v\mid \| z\|_v= T\} .
\end{equation}

\bp
\label{prohbeffwell}
\begin{enumerate}
\item Fix  a subset $S_0\subset S$ containing $S_\infty$.
 For an $|S|$-tuple $m=(m_v)$ of positive numbers,
define $$Z(m):=\prod_{v\in S_0} B_v(m_v)
\times \prod_{v\in S\setminus S_0}
S_v(m_v) .$$
Then
the family of sets $Z(m)$, $m_v\gg 1$, provided non-empty,
is effectively well-rounded.
\item The family of height balls
$B_T=\{z\in Z_S : \H_S(z)\leq T\}$ is
effectively well-rounded.
\end{enumerate}
\ep

\proof The proof relies heavily on the appendix \ref{secapp}.
For (1), we may assume that $S$ contains only one place $v$.
When $v$ is infinite, the condition
 \ref{defeffwel} \eqref{bn} is Lemma \ref{lemex6} (2)
and \ref{defeffwel} \eqref{bnhhh} is
Corollary \ref{cordevas}.c.
When $v$ is finite, the condition \ref{defeffwel} \eqref{bn} is empty
and the condition \ref{defeffwel} \eqref{bnhhh} is
Corollary \ref{corden}.b.

For (2),  \ref{defeffwel} \eqref{bn} is Lemma \ref{lemex6} (4), and
\ref{defeffwel} \eqref{bnhhh} is a combination of the following
lemma \ref{lemhhwell} with the facts that, on one hand
one has $B_{T,\e}\subset B_{kT}$
for some fixed $k>0$ and,
on the other hand,  one has $V_{kT}=O(V_T)$
again by Lemma \ref{lemex6} (4).
\hb

\begin{Lem}
\label{lemhhwell} Let
$B_T=\{z\in Z_S : \H_S(z)\leq T\}$ and $V_T:=\mu_{Z_S}(B_T)$. Then,
for any $k>0$, there exists $\delta>0$ such that
$$\int_{B_T} \H_S (z)^{-k} d\mu_{Z_S}(z) =O( V_T ^{1-\delta}) .$$
\end{Lem}

\proof
We may assume that,
for all $v$ in $S$ and $z$ in $Z_v$,
one has $\| z\|_v\geq 1$.
Set $b_T=\{ z\in Z_{S_\infty} : \H_{S_\infty}(z) \leq T\}$
and $v_T=\mu_{Z_{S_\infty}}(b_T)$.

We first claim that there exists $\delta>0$ and
$C>0$ such that for any $T>0$,
\begin{eqnarray}
\label{eqnfv2}
 \int_{b_T}
 \H_{S_\infty} ^{-k} d\mu_{Z_{S_\infty}}< C\,
v_T^{1-\delta}.
\end{eqnarray}

Set $u_T$ to be the left hand side of the above inequality.
For $T$ large, by Lemma \ref{lemex6} (1) one has
$v_T=O( T^{m_0})$ for some $m_0>0$,
hence  the derivative $u'_T=
T^{-k}v'_T$ satisfies $u'_T =O (v_T^{-k/m_0 }v'_T) $ and, integrating,
one gets $u_T=O(v_T^{1-k/m_0})$.
For $T$ bounded, since $\H_{S_\infty}$ is bounded below,
one gets $u_T=O(v_T)$.
Putting this together, one gets \eqref{eqnfv2}.

Now, for any tuple $m=(m_v)\in \m N^{S_f}$, set
 $S(m):=\prod S_v(m_v)\subset Z_{S_f}$
where
$S_v(m_v):=\{ z\in Z_v : \| z\|_v=m_v\}$.

Also set
$\pi_m:=\prod _v m_v\in \m N$ and $\om_m=\mu_{Z_{S_f}}(S(m))$.

Then we have,
where the following sums are taken over the tuples
$m\in \m N^{S_f}$ for which $S(m)$ is non empty,
\begin{eqnarray*}
\int_{B_T} \H_S (z)^{-k} d\mu_{Z_S}(z)
&=& \textstyle
\sum_{m}
\left( \int_{b_{T/\pi_m}}
 \H_{S_\infty} ^{-k} d\mu_{Z_{S_\infty}} \right)\,
\left( \int_{S(m)}
 \H_{S_f} ^{-k} d\mu_{Z_{S_f}} \right)\\
&=& \textstyle
\sum_{m} \pi_m^{-\frac{k}{2}}
\left( \int_{b_{T/\pi_m}}
 \H_{S_\infty} ^{-k} d\mu_{Z_{S_\infty}} \right)\,
\left( \int_{S(m)}
 \H_{S_f} ^{-\frac{k}{2}} d\mu_{Z_{S_f}} \right)\\
&\le & \textstyle
C\sum_{m } \pi_m^{-\frac{k}{2}}
\left(v_{T/\pi_m}\, \om_m\right)^{1-\delta}\\
 &\le & \textstyle
 C
\left(\sum_m \pi_m^{-\frac{k}{2\de}}\right)^\delta
\left(\sum_{m} v_{ T/\pi_m} \om_m\right)^{1-\delta}
\\
&\leq & \textstyle
C
\left( \prod_{v\in S_f} (1-q_v^{-\frac{k}{2\de}})\right)^{-\de}
V_T^{1-\delta}
,
\end{eqnarray*}

\noi
with positive constants $C$ and $\delta$
given by \eqref{eqnfv2} and Corollary \ref{corden}.b.
\hb

\section{Effective applications}
\label{seceffapp} In this section, assuming $K$ is a number field
and $\bh$ is a symmetric $K$-subgroup,
we give proofs of effective versions of our main theorems listed
in the introduction, keeping the notations therein.
\bc
\label{ec2} Assume that $Z_{\c O_S}\ne \emptyset$. Then
for any finite $v\in S$,
there exists $\delta>0$ such that
for any effectively well-rounded subset $\Omega
\subset Z_{S-\{v\}}$,
$$\# \{z\in Z_{\c O_S} \cap \Omega: \|z\|_v=T\} =
\wt\Vol (\Omega\times S_v(T)) \, (1+ O (T^{-\delta }))$$
as $T\to \infty$ subject to $S_v(T) \ne \emptyset$.
\ec

Recall that the normalized volume $\wt\vol$
has been defined in (\ref{eqnvolt}).
This corollary is an equidistribution statement since
one has
$$\wt\Vol (\Omega\times S_v(T))=
C\, \mu_{Z_{S-\{v\}}}\! (\Om)\,\mu_{Z_{v}}\! (S_v(T))$$
with a constant
$C$ independent of $n$ and $\Om$.\vs

\noindent{{\it Proof of Corollary \ref{ec2}}:}
The same claim with the error term $O(T^{-\delta})$ replaced by
$O((\mu_{Z_v}(S_v(T))^{-\delta})$ follows immediately from
Theorems \ref{thmeet}, \ref{thspeed} and Proposition \ref{prohbeffwell}.
Now by Corollary \ref{corden}.c, we have
a constant $a>0$
such that for $T$ large and $S_v(T)\ne \emptyset$,
 $$\mu_{Z_{v_0}}(S_{v}(T))^{-1}=O(T^{-a/2}).$$
This proves the claim.
\hb\vs

For
$z\in {\bf V}_{\mathbb Q}$, the condition $\|z\|_p\le p^n$
 is equivalent to
$ z\in p^{-n}\mathbb {\bf V}_{\z}$.
Hence we obtain:

\bc
\label{corintro2}
Assume $K=\m Q$ and fix a prime $p$
such that ${ Z}_{p}$ is non compact
and $Z_{\m Z[p^{-1}]}$ is non empty.
Then there exists $\de>0$ such that,
for any non-empty compact subset $\Omega
\subset Z_{\mathbb R}$ with piecewise smooth boundary,
$$\#\, \Omega \cap p^{-n} {\bf V}_{\z} =
\wt\vol(\Omega \times B_p(p^n))\,
(1+O(p^{-\delta n}))$$
as $n\to \infty$,
where $B_p(p^n)$ is defined in \eqref{sphered}.
\ec

Note again that
$\wt\Vol (\Omega\times B_p(p^n))=
C\, \mu_{Z_{\m R}}\! (\Om)\,\mu_{Z_{p}}\! (B_p(p^n))$
for some $C>0$ independent of $n$ and
$\Omega$.
\vs

\noindent{{\it Proof of Theorems $\ref{ec1}$ and $\ref{thoeqef}$}
Letting $v_T:=\operatorname{vol} (B_S(T))$,
Theorem \ref{ec1} with $O(T^{-\delta})$ replaced by $O(v_T^{-\delta})$
 immediately follows from Theorems \ref{thmeet} and \ref{thspeed} and
 Proposition \ref{prohbeffwell}.
To obtain the given error term, note that
at least one of the factors $Z_{v_0}$ is non-compact.
Fix $R_0>1$ such that
the volume of $B_{S\setminus \{v_0\}} (R_0)=
\{z\in Z_{S-v_0}: \H_{S\setminus \{v_0\}} (z)\le R_0\}$ is positive.
Since $B_S(T)$ contains the product $B_{v_0}(T R_0^{-1})
B_{S\setminus \{v_0\}}(R_0)$,
we have $v_T^{-1}= O( \mu_{Z_{v_0}}(B_{v_0}(T R_0^{-1}))^{-1}).$
By Corollary \ref{cordevas}.d for $v_0$ archimedean
and Corollary \ref{corden}.c for $v_0$ non-archimedean, we have
a constant $a>0$ satisfying
$$\mu_{Z_{v_0}}(B_{v_0}(T))^{-1}=O(T^{-a/2}) .$$
This proves the claim.
The same proof works for Theorem \ref{thoeqef}
applying Corollary \ref{corspeed}
in place of Theorem \ref{thspeed}.
\hb\vs

\noindent {\it Proof of Corollary $\ref{corinv1}$.}
All three cases fit in our setting as in the introduction.

For (A), if $f$ has signature $(r,s)$,
${\bf V}_1(\m R)$ can be identified with
$\op{Spin}(r-1,s)\ba \op{Spin}(r,s)$ where
the $\op{Spin}(r,s)$-action on ${\mathbb R}^{r+s}$
is given through the projection
$\op{Spin}(r,s)\to \op{SO}(r,s)$.

For (B): we have the action of $\bg={\bf SL}_n$ on
${\bf V}$ by $(g,v)\mapsto g^tvg$.
And ${\bf V}_{1}(\m R)$ is a finite disjoint union
of $\op{SO}(r, s)\ba \SL_n(\mathbb R)$
for $r+s=n$, each of them being
the variety consisting of symmetric matrices of signature $(r,s)$.

For (C), we have ${\bf V}_1(\m R)=\op{Sp}_{2n}(\mathbb R)\ba
\SL_{2n}(\mathbb R)$ with the action $(g, v)\mapsto g^tvg$.

Note for (A), if $n=3$,
$\bh_{\mathbb R}=\op{Spin} (1,1)$ may arise
and the additional assumption that
$f$ does not represent $0$ over $\mathbb Q$ implies that
$\bh$ does not allow any non-trivial $\mathbb Q$-character.
In all other cases, $\bh$ is semisimple and hence has no
character.

Now we give a uniform proof assuming that $S=\{\infty,p\}$
for the sake of simplicity. It is easy to generalize the
argument for a general $S$. Also note that this proof works equally
well for any homogeneous integral polynomial $f$ whose level
set can be identified with a symmetric variety in our set-up.

Let $d=\op{deg} f$.
For each $0\le j\le d-1$, consider
the radial projection $\pi_j: V_{p^{kd+ j}}\to  V_{p^j}$
given by $x\mapsto
p^{-k}x$.
Then
since the degree of $f$ is $d$ and the radial projection is bijective,
\begin{align*}
\#\, \Omega\cap\pi ({\bf V}_{p^{kd+j}} (\z)^{\text{prim}})
&=\# \,\{z\in {\bf V}_{p^j} (\mathbb Z [p^{-1}])\cap \Omega_j:
\|z\|_p = p^k \} \end{align*}
where $\Omega_j:=p^{j/d} \Omega\subset {\bf V}_{p^j}$.
Since ${\bf V}_{p^j}(\z [p^{-1}])$ is a finite union
of $ {\bf G} (\mathbb Z [p^{-1}]) $-orbits,
we obtain by Corollary \ref{corspeed} with
$S_0=\{\infty\}$ and $S_1=\{ p\}$
$$\#\,\Omega\cap\pi ({\bf V}_{p^{kd+j}} (\z)^{\text{prim}})\sim
\, \om_{p^{kd+j}}\vol(\Omega)   (1+ O( \om_{p^{kd+j}}^{-\delta })) $$
where $\om_{p^{kd+j}}=\mu (\{x\in {\bf V}_{p^j}
(\mathbb Q_p):\|x\|_p = p^k\})$.
Note that, by Remark \ref{remapplication},
$\om_{m}$ go to infinity with $m$
when it is non zero.
\hb\vs

\noindent{\it Proof of Corollary $\ref{corinv2}$}.
As before, we assume $S=\{\infty ,p\}$ for simplicity.
We use the same notation as in the above proof.
Then for each fixed $0\le j\le d-1$,
$f(x)=p^{kd+j}$ is equivalent to $f(p^{-k}x)=p^j$
and, if $z=p^{-k}x$ with $x\in {\bf V}(\m Z)^{\text{prim}}$, one has
$\|z\|_\infty \, \|z\|_p=\|x\|_\infty .$
Therefore
\begin{align*}N_{j, T}&
:=\#\{x\in {\bf V}(\mathbb Z)^{\text{prim}}: \|x\|_\infty<T,\;\;
f(x)= p^{kd+j}\;\;
\text{for some integer $k\geq 0$} \}\\
&=\#\{z\in {\bf V}_{p^j}(\mathbb Z[p^{-1}]): \|z\|_\infty\, \|z\|_p <T\}.
\end{align*}
By Theorem \ref{ec1}, one has
$N_{j,T}=v_{j,T}(1+O(v_{j,T}^{-\delta_j}))$
where $$v_{j,T}=\wt{\vol}(
\{(z_\infty, z_p)\in {\bf V}_{p^j}(\mathbb R)
\times {\bf V}_{p^j}(\mathbb Q_p):
 \|z_\infty\|_\infty\, \|z_p\|_p <T\}). $$

Since
$\#\{x\in {\bf V}(\mathbb Z)^{\text{prim}}: \|x\|_\infty<T, \;\;
f(x)\in p^{\m Z}\} =
\textstyle\sum_{j}N_{j, T},$
and $v_T=\sum_{j} v_{j,T}$, this proves the claim.
\hb

\section{More examples}
\label{secexample}\label{last}

Here are a few concrete examples of
applications of Theorem \ref{thoeq}
to emphasize the meaning of our results.
For each of them,
we have selected a specific global field $K$
with sets $S_0$, $S_1$
(most often $K=\m Q$,
 $S_0=\{v_0\}$ and $S_1=\{v_1\}$)
and we have selected a classical symmetric
space $Z$ defined over $K$. We look at the repartition
of $S$-integral points $z$ in $Z_{v_0}$
when imposing conditions on the $v_1$-norm of $z$.
\vs

\noi {\bf Symmetric matrices with two real places.}
This example is very classical. Let
$\tau$ be the non trivial automorphism of the real quadratic
field $K:=\m Q[\sqrt{2}]$ and set, for $d\geq 2$,
$$
Z_{\{\infty\}}:=\{ M\in \op{ M}_d(\m R)\;\;
\mbox{\rm positive definite symmetric matrix of determinant} \; 1 \}
$$
and
$$
Z_n:=\{ M\in Z_{\{\infty\}}\cap \op{ M}_d(\m Z[\sqrt{2}])\;|\;
\textstyle\sum_{i,j}\tau(m_{i,j})^2\leq n^2\}\; .
$$
\bl
\label{lemexample1}
As $n\to\infty$, these discrete sets $Z_n$
become equidistributed in
the non-compact Riemannian symmetric space $Z_{\{\infty\}}$.
\el

\proof  Let $v_0$
and $v_1$ be the two infinite places
of $K$~:
for $\la\in K$, $|\la|_{v_0}=|\la|$ and $|\la|_{v_1}=|\tau(\la)|$.
Apply Theorem \ref{thoeq} to
\addtocounter{equation}{-1}
\begin{eqnarray}
\label{eqnsetting1}
K=\m Q[\sqrt{2}]
\;,\;\;
S_0=\{ v_0\}
\; ,\;\;
S_1 =\{v_1\}
\; ,\;\;
{\bf Z}={\bf SL}_d/{\bf SO}_d
\; ,
\end{eqnarray}
and to
the group ${\bf G}={\bf SL}_d$ which acts
by $M\ra gM\,{}^tg$
on the vector space ${\bf V}$
of symmetric $d\times d$-matrices, with
${\bf Z}\sim \{  M\in {\bf V}
\;\; |\:\;
{\rm det}(M)=1
\; \}
$ as a ${\bf G}$-orbit.

Note that the group $\op{SL}(d,\m R)$ acts transitively
on $Z_{\{ \infty\} }$.
\hb
\vs

\noi {\bf Orthogonal projections with one real and one finite place.}
This example is also quite classical.
Let $p$ be a prime number, $d=d_1+d_2\geq 3$,
$$
Z_{\{\infty\}}:=\{ \pi\in \op{ M}_d(\m R)\;\;
\pi^2=^t \pi=\pi \;\; {\rm and}\;\;
{\rm tr}(\pi)=d_1 \}
$$
the Grassmannian of $\m R^d$, and
$$
Z_n:=\{ \pi\in Z_{\{\infty\}}\;|\;
p^n\pi\in \op{ M}_d(\m Z)\}\; .
$$
\bl
\label{lemexample2}
As $n\to\infty$,  these discrete sets $Z_n$
become equidistributed in
the compact Riemannian symmetric space $Z_{\{\infty\}}$.
\el

\proof  Apply Theorem \ref{thoeq} and
Remark \ref{remapplication} with
\addtocounter{equation}{-1}
\begin{eqnarray}
\label{eqnsetting2}
K=\m Q
\; ,\;\;
S_0 =\{\infty\}
\;,\;\;
S_1=\{ p\}
\; ,\;\;
{\bf Z}={\bf O}_d/{\bf O}_{d_1}\times {\bf O}_{d_2}
\; ,
\end{eqnarray}
and to
the group ${\bf G}={\bf Spin}_d$ which acts,
by conjugation via ${\bf SO}_d$,
on the vector space ${\bf V}$
of $d\times d$-matrices, with
${\bf Z}\sim \{  \pi\in {\bf V}
\;\; |\:\;
\pi^2=^t \pi=\pi \;\; {\rm and}\;\;
{\rm tr}(\pi)=d_1 \}
$ as a ${\bf G}$-orbit.

Note that the group $\op{Spin}(d,\m R)$ acts transitively
on $Z_{\{ \infty\} }$.
\hb
\vs

\noi {\bf Complex structures with one finite and one real place.}
In this example, one chooses a prime number $p$
and set, for $d\geq 1$,
$$
Z_{\{ p\}}:=\{ J\in \op{ M}_{2d}(\m Q_p)
\;\; |\:\;
J^2=-{\rm Id}
\;\; {\rm and}\;\; tr(J)=0
\; \}
$$
and
$$
Z_R:=\{ J\in Z_{\{ p\}}\cap\op{ M}_{2d}(\m Z[\textstyle\frac1p])
\;\; |\;\;
\sum_{i,j}J_{i,j}^2\leq R^2\}\; .
$$
\bl
\label{lemexample3}
As $R\to\infty$, these discrete sets $Z_R$
become equidistributed in
the $p$-adic symmetric space $Z_{\{ p\}}$.
\el

\proof  Apply Theorem \ref{thoeq} to
\addtocounter{equation}{-1}
\begin{eqnarray}
\label{eqnsetting3}
K=\m Q
\; ,\;\;
S_0 =\{ p\}
\;,\;\;
S_1=\{ \infty\}
\; ,\;\;
{\bf Z}={\bf GL}_{2d}/{\bf GL}_d\times {\bf GL}_d
\; ,
\end{eqnarray}
and to
the group ${\bf G}={\bf SL}_{2d}$ which acts
by conjugation on the vector space ${\bf V}$
of $2d\times 2d$-matrices, with
${\bf Z}\sim \{  J\in {\bf V}
\;\; |\:\;
J^2=-{\rm Id}
\;\; {\rm and}\;\; tr(J)=0
\; \}
$ as a ${\bf G}$-orbit.

Note that the group $\op{SL}(2d,\m Q_p)$ acts transitively
on $Z_{\{ p\} }$.
\hb
\vs

\noi {\bf Antisymmetric matrices with two finite places
in characteristic zero.}
In this example, one chooses two distinct prime numbers
$p$ and $\ell$,
and set, for $d\geq 1$, $n\geq 0$, and $R>d$,
$$
Z_{\{ p\}}:=\{ A\in \c M_{2d}(\m Q_p)
\; |\;
A=-^t A
\;\; {\rm and}\; \;
{\rm det}(A)=1\; \}
$$
and
$$
Z_{n,R}:=\{ A\in Z_{\{ p\}}
\;|\;
\textstyle
\sum_{i,j}A_{i,j}^2 < R
\;\; {\rm and}\;\;
\ell^nA\in \op{ M}_{2d}(\m Z[\textstyle\frac{1}{p}])\}\; .
$$
\bl
\label{lemexample4}
As $n+R\to\infty$, these discrete sets $Z_{n,R}$
become equidistributed in
the $p$-adic symmetric space $Z_{\{ p\}}$.
\el

\proof  Apply Theorem \ref{thoeq} and
Remark \ref{remapplication} with
\addtocounter{equation}{-1}
\begin{eqnarray}
\label{eqnsetting4}
K=\m Q
\;,\;\;
S_0 =\{p\}
\; ,\;\;
S_1=\{ \infty, \ell \}
\; ,\;\;
{\bf Z}={\bf SL}_{2d}/{\bf Sp}_{d}
\; ,
\end{eqnarray}
and to
the group ${\bf G}={\bf SL}_{2d}$ which acts
by $g\ra g A\, {}^tg$ on the vector space ${\bf V}$
of antisymmetric $2d\times 2d$-matrices, with
${\bf Z}\sim \{  A\in {\bf V}
\;\; |\:\;
{\rm det}(A)=1\; \}
$ as a ${\bf G}$-orbit.

Note that the group $\op{SL}_{2d}(\m Q_p)$ acts transitively
on $Z_{\{ p\} }$.\hb
\vs

\noi {\bf Quadrics with two places
in positive characteristic.}
In this example, $p$ is an odd prime, and
one set, for $d\geq 3$,
$$
Z_{\{ 0\}}:=\{ P\in \m F_p((t))^d\;\;|\;\;
P_1^2+\cdots +P_d^2=1\; \}
$$
and, for $n_1\geq n_2\geq \ldots \geq n_d\geq 0$,
$$
Z_{n_1,\ldots , n_d}
:=\{ P\in Z_{\{ 0\}}\cap \m F_p[t,t^{-1}]^d
\;\;|\;\;
 {\rm deg}(P_i)=n_i
\;\;\forall i
\;
\}\; ,
$$
where deg$(\sum a_it^i):={\rm max}\{ i\;|\; a_i\neq 0\}$.

\bl
\label{lemexample5}
If   $n_1=n_2=n_3$ goes to infinity
or if $p\equiv 1 \; {\rm mod}\; 4$ and
$n_1=n_2$ goes to infinity, these discrete sets $Z_n$
become equidistributed in
the sphere $Z_{\{ 0\}}$.
\el

\proof
Let $0$ and $\infty$
be the two (finite) places of the field $\m F_p(t)$
associated to the two points $0$ and $\infty$
of $\m P^1(\m F_p)$.
Apply Theorem \ref{thoeq} and
Remark \ref{remapplication}
with
\addtocounter{equation}{-1}
\begin{eqnarray}
\label{eqnsetting5}
K=\m F_p(t)
\; ,\;\;
S_0 =\{ 0 \}
\;,\;\;
S_1=\{ \infty\}
\; ,\;\;
{\bf Z}={\bf SO}_{d}/{\bf SO}_{d-1}
\; ,
\end{eqnarray}
and to
the group ${\bf G}={\bf Spin}_{d}$ which acts naturally,
via ${\bf SO}_{d}$,
on the $d$-dimensional vector space ${\bf V}$, with
the sphere
${\bf Z}\sim \{  v\in {\bf V}
\;\; |\:\;
v_1^2+\ldots v_d^2=1\; \}
$ as a ${\bf G}$-orbit.

The corresponding two completions
are the fields of Laurent series
$K_{0}=\m F_p((t))$ and $K_{\infty} =\m F_p((t^{-1}))$,
and the ring of $S$-integers is $\c O_S=\m F_p[t,t^{-1}]$.
Set
$$
Z_{\{\infty\}}:=\{ P\in \m F_p((t^{-1}))^d\;\;|\;\;
P_1^2+\cdots +P_d^2=1\; \} \; ,
$$
note  that the well-rounded subset
$$
B_{n_1,\ldots , n_d}:= \{ P\in Z_{\{ \infty\}}
\;\;|\;\;
 {\rm deg}(P_i)=n_i
\;\;\forall i
\;
\}
$$
is non-empty if and only if
$n_1=n_2=n_3$
or $p\equiv 1 \; {\rm mod}\; 4$ and
$n_1=n_2$.

Note also that the group ${\rm Spin}(d,\m F_p((t)))$
acts transitively
on $Z_{\{ 0\} }$.
\hb
\vs

\noi {\bf Other examples.}
The reader may construct easily many similar examples
choosing other triples $(K,S,Z)$.
For instance, ``Quadrics with three infinite places'',
``Lagrangian decompositions with two infinite and three finite place'',
``Hermitian matrices with
four places
in positive characteristic'',
and so on....

\section{Appendix: Volume of balls}
\label{secapp}
In this appendix we prove precise estimates for the volume of balls
which are needed in sections \ref{secwellrounded} and \ref{seceffwel}.
These estimates will be consequences of the following two
general theorems \ref{thje} and \ref{thden}.

\subsubsection*{\bf Volume of balls over the reals}
We will first need a variation
of a theorem on fiber integration.
This theorem says that the volume of the fibers of an analytic
function has a, term-by-term differentiable
asymptotic expansion in the scale of functions
$t^j (\log t)^{k}$ with $j$ rational and $k$ non-negative integer.
More precisely,

\bt
\label{thje}
\cite{Je}
\label{thdevas}
Let $X\subset \m R^m$ be a smooth real analytic variety,
$f:X\ra \m R$ a real analytic function
and $\nu$ a $C^\infty$ measure on $X$.
Then, for any compact $K$ of $X$,
there exist $d\in \m N$ and
a set
 $\{A_{j,k}: j\in \frac{1}{d} \m N, k\in \m Z, 0\le k< m\}$ of
distributions on $X$ supported by $f^{-1}(0)$
such that, for every $C^\infty$ function
$\ph:X\ra \m R$ with support in $K$, the integral
$$v_\ph(t):=\int_{0\leq f(x)\leq t}\ph (x)
\, d\nu(x)$$
has a term-by-term differentiable asymptotic expansion
when $t>0$ goes to $0$
$$\sum_{{j\in \frac{1}{d}\m N}}
\;
 \sum_{0\leq k< m}A_{j,k}(\ph)\, t^{j}(\log t)^{k}$$
\et

This means that, for every $j_0\geq 0$,
defining $v_{\ph,j_0}$ by truncating the
above sum:
$$v_{\ph,j_0}(t):=
\sum_{\stackrel{\mbox{\tiny $0\! <\! j\leq j_0$}}{j\in \frac{1}{d}\m N}}
\;
\sum_{0\leq k< m}A_{j,k}(\ph)\, t^{j}(\log t)^{k},$$
 one has
$\left(\frac{d}{dt}\right)^\ell (v_\ph-v_{\ph,j_0})=o(t^{j_0-\ell})$
for every $\ell\geq 0$.
\vs

\noi{\bf Remark}
This theorem is stated for a smooth analytic variety
$Z$ and a smooth measure $\mu$.
Its proof is based on the real version of Hironaka's
resolution of singularity
as in \cite{At}.
Using once more Hironaka's theorem it can be applied
to a singular analytic variety $X$ with
a measure $\mu$  associated to a meromorphic
differential form. Here is one example of such an application:

\bp
\label{prodevas}
Let $Z$ be a connected component of the real points
of a smooth real affine algebraic variety $\bf Z$
and $\mu$ a measure on $Z$ which is defined by a regular
differential form of ${\bf Z}$.
Let $F:Z\ra \m R$ be a positive proper regular function
and set $v_T:=\mu (\{ z\in Z\; :\; F(z)\leq T\})$.
Then,
there exist positive integers $\ell_0$, $d$ and constants $a_{j,k}$
such that $v_T$
has a term-by-term differentiable asymptotic expansion
when $T\ra \infty$
$$\sum_{\stackrel{\mbox{\tiny $j\leq \ell_0$}}{j\in \frac{1}{d}\m Z}}
\sum_{0\leq k< m}a_{j,k} T^{j}(\log T)^{k}.$$
\ep

The condition $F$ {\it regular} means that $F$
is the restriction to $Z$ of a regular function
on the algebraic variety ${\bf Z}$.
\vs

\proof
Using the resolution of singularities
we can view $\bf Z$ as an open real algebraic subvariety
of a smooth projective variety $\bf X$ such that
the boundary ${\bf Y}:={\bf X}-{\bf Z}$ is a divisor with normal crossing.
Hence, in a neighborhood $U_{y_0}$
of each real point $y_0$ of ${\bf Y}$,
there is a real local coordinate system $(x_1,\ldots ,x_m)$
such that ${\bf Y}$ is given by $x_1\cdots x_r=0$, for some $r\leq m$.
We are only interested in those points $y_0$ in the closure of $Z$.
Near these points, the meromorphic function $f:=1/F$ is zero
on ${\bf Y}$.
Using a partition of unity associated to such a cover,
we are reduced to a local problem.
Namely proving, for  every $C^\infty$ function $\ph$ with compact
support in $U_{y_0}$, the existence of a term-by-term
differentiable asymptotic expansion  for
$$v_{\ph,\mu}(t):=\int_{f(x)\leq t}\ph (x)
\, d\mu(x)$$
when $t:=1/T$ goes to $0$.
It is equivalent to prove the existence of such an asymptotic
expansion for the derivative
$v'_{\ph,\mu}(t)$ which is called the
{\it integral of $\mu$ on the fiber
$f^{-1}(t)$}.

Using once more the resolution of singularities
for the numerator and denominator of $f$ and a new partition of unity,
we can assume that
$f$ is monomial
in these coordinate systems.
Using the fact that $f$ is positive on $Z$ and zero on its boundary,
we deduce that
$f$ is given by $f=s x_1^{p_1}\cdots x_r^{p_r}$ with
${p_1},\ldots ,p_{r}$ positive integers and $s=\pm 1$.
Hence $f$ is an analytic function near $y_0$.

The integral  $v_{\ph,\mu}$ is now very similar to the integral
$v_{\ph}$ of Theorem \ref{thdevas}
except that $\mu$ may not be smooth.
However $\mu$ is defined by a regular differential form on $\bf Z$
hence there exists a positive integer $\ell_0$ such that
the measure $\nu:=f^{\ell_0}\mu$ is smooth.
According to Theorem \ref{thdevas},
$v_{\ph,\nu}$ has a term-by-term
differentiable asymptotic expansion.
The following equality between the {\it fiber integrals}
$$
v'_{\ph,\mu}(t)= t^{-\ell_0}v'_{\ph,\nu}(t)
$$
implies that
$v_{\ph,\mu}$ has also a term-by-term
differentiable asymptotic expansion.
\hb\vs

For us, the main example to which we will apply
Proposition \ref{prodevas} is
a closed orbit $Z$ under the group of
$\m R$-points of a $\m R$-algebraic
group,
an invariant measure $\mu$ on this orbit
and the restriction $F$ to $Z$ of the square
of an euclidean norm on $\m R^m$.
Hence we get,

\bc
\label{cordevas}
Let $Z$ be a closed orbit of the group $G$ of $\m R$-points
of an $\m R$-algebraic group
acting algebraically on a $\m R$-vector space  $V$,
$\mu$ an invariant measure on $Z$ and
$\|\cdot \|$ an euclidean norm on $V$.
Set $B_T:=\{z\in Z: \|z\| \leq T\}$ and $v_T:=\mu( B_T)$.\\
a) Then
$v_T
\sim c \, T^a (\log T)^{b}$, as $T\to \infty$,
where $a\in \q_{\ge 0}$, $b\in \m Z_{\ge 0}$  and $c>0$.\\
b) Moreover
$\frac{d}{dT}v_T
\sim c \frac{d}{dT}( T^{a} (\log T)^{b})$, as $T\to \infty$.\\
c) For any $k_0>0$, there exists $\delta_0>0$ such that one has,
as $T\to \infty$,
$$
\int_{B_T}\|z\|^{-k_0} d\mu (z) =
O(v_T^{1-\delta_0}) .
$$
d) If $G$ is semisimple and $Z$ is non compact then one has $a\neq 0$.
\ec

\noi {\bf Remarks } - When $Z$ is a symmetric variety,
the point a) is proven in \cite[Corollary 6.10]{GOS1} for any norm on $V$
and
the parameters $a$ and $b$ are
explicitly given.\\
- When $G$ is a group of diagonal matrices, the constant $a$ is zero.
\vs

\proof
a) and b) This is a special case of
Proposition \ref{prodevas}. Note that since $v_T$ is an increasing function
of $T$, one has $a\geq 0$.
Moreover, note that, when $a=b=0$, the orbit is of finite volume
hence compact.

c) Set $u_T:=\int_{B_T}\|z\|^{-k_0} d\mu (z)$.
By a), one has
$v_T=O( T^{m_0})$ for some $m_0>0$.
Hence  the derivative $u'_T=
T^{-k_0}v'_T$ satisfies $u'_T =O (v_T^{-k_0/m_0 }v'_T) $ and, integrating,
one gets $u_T=O(v_T^{1-k_0/m_0})$.

d) This is a special case of the following Proposition \ref{proazero}
\hb

\bp
\label{proazero}
With the notations of Corollary \ref{cordevas},
one has the equivalence:
\centerline{
All unipotent elements of $G$ act trivially on
$Z$ $\Longleftrightarrow a=0  $.}
\ep

\proof $\Longrightarrow$
By assumption the normal subgroup of $G$ generated
by the unipotent elements of $G$
acts trivially on $Z$.
Hence one can assume that $G$ is a product of a compact group
by a $r$-dimensional group of diagonal matrices.
In this case one has $\mu (B_T)= O((\log T)^r)$ as $T\to\infty$.

$\Longleftarrow$
This implication is a consequence of the following Lemma \ref{lemazero}.
\hb\vs

\bl
\label{lemazero}
Let $U$ be a one-parameter unipotent subgroup of ${\rm GL}_m(\m R)$,
$\mu$ a $U$-invariant measure on $\m R^m$
which is not supported by the $U$-fixed points
and set $B_T$ for the euclidean ball of radius $T$ on $\m R^m$.
Then one has
$$\D\liminf_{T\ra\infty} \frac{\log(\mu (B_T))}{\log(T)} >0.$$
\el

\proof First of all, note that all the orbits $Uz$ of $U$
in $\m R^m$ are images of $\m R$ by polynomial maps $t\mapsto u_t z$
of degree $d_z\leq m$.
We may assume that this degree $d_z$ is $\mu$-almost everywhere
non-zero constant.
Set $d\geq 1$ for this degree, write
$u_t z= t^d v_z + O(t^{d-1})$ for some non-zero $v_z\in \m R^m$,
and note that the constant involved
in this $O(t^{d-1})$ is uniform on compact
subsets of $\m R^m$.

One can find a compact subset $C\subset\m R^m$
transversal to the $U$-action such that
$\mu(UC)>0$.
The pull-back on $\m R\times C$ of the measure $\mu$
by the action
$ (t,z)\mapsto  u_t z$ has the form $dt \otimes \nu$ where $dt$
is the Lebesgue measure on $\m R$ and $\nu$ is a non-zero measure on $C$.
Choose $c>\sup_{z\in C}\| v_z\|$.
Then, for $R$ large, one has
$$u_{[0,R]}( C)\subset B_{cR^d}$$ and hence
$\mu(B_{cR^d})\geq R\,\nu(C)$.
This proves  our claim.
\hb

\subsubsection*{\bf Volume of balls over the $p$-adics}

We will also need Denef's theorem on $p$-adic integration.
For that we need some notations.
A subset of $\m Q_p^m$ is said {\it semialgebraic}
if it is obtained by boolean
operations from sets $P_{f,r}:=\{ x\in \m Q_p^m\;
/\;\exists y\in \m Q_p\; :\;
f(x)=y^r\}$ with $f$ a polynomial in $m$ variables
with coefficients in $\m Q_p$ and $r\geq 2$.
According to Macintyre's theorem, which is the $p$-adic analog of
Tarski-Seidenberg theorem, those sets are exactly the definable
sets of the field $\m Q_p$ \cite{Mac}.
A function $f$ between two $\m Q_p$-vector spaces
is said {\it semialgebraic}
if its graph is semialgebraic.
According to Denef's cell decomposition theorem (\cite{De2} and \cite{Clu}),
for every semialgebraic
subset $S$, there exists a finite partition of $S$
in semialgebraic sets $S_1,\ldots ,S_{j_{max}}$
(called cells) such that, for each $j=1,\ldots ,j_{max}$,
$S_j$ is in semialgebraic bijection
with a semialgebraic open subset $O_j$ of a
vector space $\m Q_p^{d_j}$
(recently, R. Cluckers has shown the
existence of a semialgebraic bijection between $S$ itself
and some $\m Q_p^d$).
A measure $\mu$ on $S$ is said  {\it semialgebraic}
if there exists a cell decomposition of $S$
on each cell of which $\mu$ is of the form $|g_j(x)|dx$
where $g_j$ is a semialgebraic function
on $\m Q_p^{d_j}$ and $dx$ is a Haar measure on $\m Q_p^{d_j}$.
A function $a:\m Z\ra\m Z$ is said {\it simple}
if there are finite partition of $\m N$ and $-\m N$
by finite sets and arithmetic progressions on which $a$ is affine,
see \cite[\S 2.13, 2.14 and 4.4]{De1}.

\bt
\cite[Theorem 3.1]{De1}
\label{thden}
Let $\mu$ be a semialgebraic measure
on an $m$-dimensional semialgebraic subset $S$ over $\m Q_p$
and $f$ be a semialgebraic
function on $S$. For $n\in \m Z$, set
$$I_n:=\int_{|f(x)|=p^n}\,d\mu(x)$$
when this integral is finite and $I_n=0$ otherwise.
Then, for all $n\in \m Z$, one has
$$I_n=\sum_{1\leq i\leq e}\ga_i(n)p^{\be_i(n)}$$
where $e\in \m N$, $\be_i:\m Z \to \m Z$ is a simple function and
$\ga_i:\m Z \to \m Z$ is a product of at most
$m$ simple functions for each $1\le i\le e$.
\et
For instance,
an orbit under the group of $\m Q_p$-points of a $\m Q_p$-algebraic
group acting algebraically is definable and hence  semialgebraic,
by Macintyre's theorem,  and
an invariant measure on this orbit is semialgebraic.
Hence one gets:

\begin{Cor}
\label{corden}
Let $k$ be a finite extension of $\m Q_p$,
$q$ the absolute value of an uniformizer,
$G$ the group of $k$-points
of an algebraic $k$-group,
$\rho:G\ra GL(V)$ a representation of $G$ defined over $k$,
$Z$ a closed $G$-orbit in $V$,
$\mu$ an invariant measure on $Z$ and
$\|\cdot \|$ a max norm on $V$.
Denote by $S_T$ the sphere
$ S_T=\{z\in {Z}: \|z\|=T\} $ and set $v_T:=\mu(S_T)$.\\
a) There exists $N_0\in \n$ such that, for
each $0\le j_0 < N_0$, one of the following holds:
\begin{enumerate}
\item $ S_{q^j}$ is empty, for $j\equiv j_0 \mod N_0$ large;
\item there exist $a_{j_{_0}}\in \q_{\ge 0}$, $b_{j_{_0}}\in \m Z_{\ge 0}$,
and $c_{j_{_0}}>0$
such that,
$$v_{q^j}\sim c_{j_{_0}}\, q^{\, a_{j_{_0}} j} \,j^{\, b_{j_{_0}}}
\;\;\;\;\;\; \;\;\;\;\;\; {\rm for}\;\;
j\equiv {j_{_0}} \mod N_0\;\; \; large.$$
\end{enumerate}
b) For any $k_0>0$, there exists $\delta_0>0$ such that one has,
as $T\to\infty$,
$$
\int_{S_T}\|z\|^{-k_0} d\mu (z) =
O(v_T^{1-\delta_0}).
$$
c) If $G$ is semisimple and $Z$ is non compact then,
for  all  $j_{_0}$ in case $(2)$, one has
$a_{j_{_0}} \neq 0$.
\end{Cor}

\noi {\bf Remarks }
- Let us recall that a max norm is a norm
given in some basis $e_1,\ldots,e_m$ by
$\|\sum x_ie_i\|=\max |x_i|$.\\
- When $G$ is a group of diagonal matrices, all the constants
$a_{j_{_0}}$ are zero.
\vs

\proof Viewing $V$ as a $\m Q_p$ vector space, we may
assume that $k=\m Q_p$.

a) This is a special case of
Theorem \ref{thden}.

b) By a), there exists $m_0>0$ such that
$v_T=O( T^{m_0})$.
Hence one has\\ $\int_{S_T}\|z\|^{-k_0} d\mu (z)=
T^{-k_0}v_T =O (v_T^{1- k_0/m_0 }) .$

c) This is a special case of the following Proposition \ref{proazero2}
which is analogous to Proposition \ref{proazero}.
\hb

\bp
\label{proazero2}
With the notations of Corollary \ref{corden},
the following are equivalent:

$(i)$ All unipotent elements of $G$ act trivially on $Z$,

$(ii)$ For all $j_{_0}$ in case $(2)$, one has $a_{j_{_0}}=0$,

$(iii)$ Either $Z$ is compact or, for some $j_{_0}$
in case $(2)$, one has $a_{j_{_0}}=0$.
\ep

\proof The proof is as in Proposition \ref{proazero}.

$(i)\Rightarrow (ii)$
By assumption the normal subgroup of $G$ generated
by the unipotent elements of $G$
acts trivially on $Z$.
Hence one can assume that $G$ is a product of a compact group
by an $r$-dimensional group of diagonal matrices.
In this case, one has $\mu(S_{p^j})=O(j^{\,r})$ as $j\to\infty$.

$(ii)\Rightarrow (iii)$ If $Z$ is non compact, at least one $j_{_0}$
is in case $(2)$.

$(iii)\Rightarrow (i)$
This implication is a consequence of the following Lemma \ref{lemazero2}.
\hb

\bl
\label{lemazero2}
Let $k$ be a finite extension of $\m Q_p$,
$U$  a one-parameter unipotent subgroup of ${\rm GL}(m,k)$,
$\mu$ a $U$-invariant measure on $k^m$
which is not supported by the $U$-fixed points
and denote by $S_T$  the sphere of radius $T$ on $k^m$
for the max norm.
Then one has, as $T\to\infty$ subject to the condition
$\mu(S_T)\neq 0$,
$${\rm liminf} \frac{\log(\mu (S_T))}{\log(T)} >0.$$
\el

\proof The proof is as in Lemma \ref{lemazero}.
First of all, note that
all the orbits $Uz$ of $U$
in $k^m$ are images of $k$ by polynomial maps
$t\mapsto u_tz$ of degree $d_z\leq m$.
We may assume that this degree $d_z$ is $\mu$-almost everywhere
non-zero constant.
Set $d\geq 1$ for this degree and write
$u_t z= t^d v_z + O(t^{d-1})$ for some non-zero $v_z\in k^m$.
Let
$q$ be the absolute value of an uniformizer.
The set
$\{ j\in \m N \; :\; Uz\cap S_{q^j}\neq \emptyset\}$ is then
equal, up to  finite sets, to some arithmetic progression $j_z+\m N d$
with $0\leq j_z< d$.
We may assume that this integer $j_z$ is $\mu$-almost everywhere
constant.
Set $j_{_0}$ for this integer.

One can find a compact subset $C\subset k^m$
transversal to the $U$-action such that
$\mu(UC)>0$ and on which $\|v_z\|$ is constant equal
to some power $q^{j_{_0}+dm_{_0}}$ with $m_{_0}\in \m N$.
The pull-back on $k\times C$ of the measure $\mu$
by the action
$ (t,z)\mapsto  u_t z$ has the form $dt \otimes \nu$ where $dt$
is a Haar measure on $k$ and $\nu$ is a non-zero measure on $C$.
Then, for $|t|=q^\ell$ large, one has
$$u_{t}( C)\subset S_{q^{j_{_0}+dm_{_0}+d\ell}}$$ and hence
$\mu(S_{q^{j_{_0}+dm_{_0}+d\ell}})\geq (q\! -\! 1)\, q^{\ell-1}\,\nu(C)$.
This proves  our claim.
\hb

\newpage


\end{document}